\documentclass{amsart}

\usepackage{graphicx}
\usepackage{amscd}
\usepackage{latexsym}
\usepackage{amssymb}
\usepackage{amsthm}
\usepackage[mathscr]{eucal}
\usepackage{amsmath}

\numberwithin{equation}{section}
\numberwithin{figure}{section}

\theoremstyle{plain}
\newtheorem{theorem}{Theorem}[section]
\newtheorem{corollary}[theorem]{Corollary}
\newtheorem{lemma}[theorem]{Lemma}
\newtheorem{proposition}[theorem]{Proposition}

\theoremstyle{definition}
\newtheorem{definition}[theorem]{Definition}

\theoremstyle{remark}
\newtheorem{remark}[theorem]{Remark}
\newtheorem{example}[theorem]{Example}

\begin{document}

\title{The t-singular homology of orbifolds}

\author{Yoshihiro Takeuchi}
\address{Department of Mathematics,
         Aichi University of Education,
         Igaya, Kariya 448-0001, Japan}
\email{yotake@auecc.aichi-edu.ac.jp}

\author{Misako Yokoyama}
\address{Department of Mathematics,
         Faculty of Science, Shizuoka University,
         Ohya, Shizuoka 422-8529, Japan}
\email{smmyoko@ipc.shizuoka.ac.jp}


\begin{abstract}
For an orbifold $M$ we define a homology group, called t-singular homology group $t$-$H_q(M)$, which depends not only on the topological structure of the underlying space of $M$, but also on the orbifold structure of $M$. We prove that it is a b-homotopy invariant of orbifolds. If $M$ is a manifold, $t$-$H_q(M)$ coincides with the usual singular homology group. We calculate the t-singular homology groups with $\mathbb{Z}$-coefficients of several orbifolds. The t-singular homology group with rational coefficients of an orbifold $M$ is isomorphic to the singular homology group with rational coefficients of $M$ (and to that of $|M|$).

For each $q$ we define a homomorphism from the orbifold homotopy group $\pi_q(M,x_0)$ to $t$-$H_q(M)$ by using of t-modifications. Especially, $t$-$H_1(M)$ is isomorphic to the abelianization of $\pi_1(M,x_0)$ if $M$ is arcwise connected.
\end{abstract}

\maketitle

\section{Introduction}

The notion of orbifold was first introduced by [Sa] under the name of $V$-manifold. The name of \lq orbifold\rq\; is after [Th], where orbifolds, especially 3-orbifolds, were studied hardly from a different view point.

Speaking of the algebraic invariants of orbifolds, in [Sa] Satake studied the singular homology, de Rham cohomology and check (co)homology. The singular homology group of an orbifold is isomorphic to the singular homology group of its underlying space. In [Th] Thurston studied the Euler number, and the fundamental group of an orbifold. The latter was defined as the deck transformation group of the universal orbifold covering space. In \cite{F-S} Furuta and Steer defined a homotopy group and a cohomology of an orbifold.

In this paper we define a homology group of an orbifold $M$, called t-singular homology group $t$-$H_q(M)$, which respects the orbifold structure and is not a topological invariant of the underlying space in general. The homotopy group $\pi_q(M,x_0)$ of $M$ is regarded as the group of homotopy classes of $q$-dimensional singular spheres with base point $x_0$. Then we define the Hurewicz homomorphisms from the orbifold homotopy group to the t-singular homology group.

When we imagine what \lq a singular homology group of an orbifold $M$' is, three types of those might come into our minds. The first one is the singular homology group of the underlying space of $M$. The second one is the singular homology group of a topological space constructed from the orbifold $M$, or its variation. The third one is the singular homology group obtained by a generalization of the singular homology theory of manifolds to that of orbifolds. Here in this paper we will choose the third one.

First we prepare two kinds of continuous maps of orbifolds. One is called just a continuous map, which consists of a continuous map between the underlying spaces and a family of continuous maps between local uniformizations. The other is called a b-continuous map, which is continuous and maps at least one regular point to a regular point. The latter map was called an orbi-map in \cite{japanfinite,finite,least}, etc. Then homotopy equivalence and b-homotopy equivalence are considered. A b-homotopy was called an orbi-homotopy in \cite{japanfinite,finite,least}, etc. Isomorphisms and covering maps between orbifolds are regarded as special cases of continuous maps of orbifolds.

Next we define a $q$-dimensional t-singular simplex of an orbifold $M$ to be a continouous map from the standard $q$-dimensional simplex to $M$, which is \lq tame\rq\; and \lq transverse\rq\; with respect to the singular set of $M$. Then the t-singular homology group of $M$ is defined in a usual way by using of the chain complex of the free abelian group generated by t-singular simplices of $M$. Since an isomorphism between orbifolds is a chain map and induces an isomorphism between their t-singular chain complices, this homology group is an isomorphism invariant of orbifolds.

In general, a b-continuous map does not induce a homomorphism between the t-singular chain complices because the composition of a t-singular simplex and a b-continuous map is not necessarily a t-singular simplex. But we can define an induced homomorphism between the t-singular homology groups by a t-modification. The homomorphism does not depend on the choice of a t-homotopy in the modification. Then we prove the b-homotopy invariance of the t-singular homology group.

The fundamental group of an orbifold is regarded as \lq the group of the b-homotopy classes of loops with fixed base point\rq. It is ismorphic to the group defind by Thurston in \cite{Th} as the fundamental group of an orbifold, that is, the deck transformation group of the universal covering space. It is also isomorphic to the group of equivalent classes of paths in the universal covering space, see \cite{finite}. We show \lq the Hurewicz isomorphism theorem' which is the abelianization from the fundamental group onto the first t-singular homology group.

In \cite{F-S} Furuta and Steer defined $\pi_q$ of an orbifold as the $\pi_{q-1}$ of the orbifold loop space. In this paper $\pi_q$ is regarded as the group of the homotopy classes of singular $q$-spheres, which is eqivalent to the above. If $q\geq 2$, then $\pi_q$ is commutative. We define the Hurewicz homomorphism from $\pi_q$ to the $q$-th t-singular homology group $t$-$H_q(M)$.

In the following we summarize the contents of the paper. In Section 2 we give some basic preparations on orbifolds. In Section 3 we describe on (b-)continuous maps, (b-)homotopies, and the fundamental groups of orbifolds. In Section 4 we define the t-singular homology group of an orbifold. In Section 5 we define a t-modification to construct a homomorphism between the t-singular homology groups induced by a b-continuous map (see Theorem \ref{f*def-th}). Theorem \ref{b-homotopy-inv-th} says that if two orbifolds are b-homotopy equivalent, then their t-singular homology groups are isomorphic. In Section 6 we define the Hurewicz homomorphism from $\pi_q(M,x_0)$ to $t$-$H_q(M)$ where $\pi_q(M,x_0)$ is the $q$-th homotopy group of an orbifold $M$ with base point $x_0$. We show the Mayer-Vietoris exact sequences and the K\"unneth's formula for the t-singular homology in Sections 7 and 8. Some examples of the t-singular homology groups are given in Section 9. The singular (not t-singular) homology group of an orbifold is isomorphic to the singular homology group of its underlying space. We point out this fact in Section 10. In Section 11 we show that the t-singular homology group with rational coefficients of an orbifold $M$ is isomorphic to the singular homology group with rational coefficients of $M$ (and to that of $|M|$). In Section 12, we give the definition of the ws-singular cohomology of an orbifold. In \cite{ws-coh} we will prove the duality theorem of the t-singular homology group between the ws-singular cohomology group.

\vskip5mm
\noindent
{\bf Acknowledgements.}\quad
The authors would like to thank to Toshiyuki Akita, \\
\noindent
Yoshihiro Fukumoto, Mikio Furuta, Mitsuyoshi Kato, Norihiko Minami, Masaharu Morimoto, and Akihiro Tsuchiya for their useful comments.

\section{Preliminaries on orbifolds}

By an {\it n-dimensional\/} ({\it topological\/}) {\it orbifold $M$}, we shall mean a Hausdorff space $X$ together with a system $\mathcal{S}=(\{U_i\},\{\varphi_i\},\{\tilde{U}_i\},\{G_i\},\{\tilde{\varphi}_{ij}\},\{\eta_{ij}\})$ which satisfies the following:
\begin{enumerate}
\item[(i)] $\{ U_i\}$ is an open cover of $X$, which is locally finite and closed under finite intersections.
\item[(ii)] For each $U_i$, there exist a finite group  $G_i$ acting smoothly and effectively on a connected open subset $\tilde{U}_i$ of $\mathbb{R}^n_+$ and a homeomorphism $\varphi_i:\tilde{U}_i/G_i\rightarrow U_i$.
\item[(iii)] If $U_i\subset U_j$, then there exists a monomorphism $\eta_{ij}:G_i\rightarrow G_j$ and an embedding $\tilde{\varphi }_{ij}:\tilde{U}_i\rightarrow \tilde{U}_j$ such that for any $\sigma \in G_i$ and $x\in \tilde{U}_i$, $\tilde{\varphi }_{ij}(\sigma x)=\eta_{ij}(\sigma )\tilde{\varphi }_{ij}(x)$ and the following diagram commutes, where $\varphi_{ij}$ is induced by the monomorphism $\eta_{ij}$ and the embedding $\tilde{\varphi}_{ij}$, and $r_i$, $r_j$ are the natural projections.

\begin{equation}
\CD
\tilde{U}_i @>\tilde{\varphi}_{ij}>> \tilde{U}_j \\
@Vr_iVV @VVr_jV \\
\tilde{U}_i/G_i @>\varphi_{ij}>> \tilde{U}_j/G_j \\
@V\varphi_iVV @VV\varphi_jV \\
U_i @>>> U_j
\endCD
\end{equation}

\item[(iv)] If $U_i\subset U_j\subset U_k$, then $\tilde{\varphi}_{jk}\circ\tilde{\varphi}_{ij}=\tilde{\varphi}_{ik}$. (From this formula and (iii), it holds that $\eta_{jk}\circ\eta_{ij}=\eta_{ik}$.)
\end{enumerate}

Each $\varphi_i\circ r_i:\tilde{U}_i\rightarrow U_i$ is called a {\it local chart} of $M$, and $\{\varphi_i\circ r_i:\tilde{U}_i\rightarrow U_i\}$ is called a {\it system of local charts} of $M$. $\mathcal{S}$ is called an {\it atlas} of $M$. If $\mathcal{S}'$ is obtained from $\mathcal{S}$ by changing each $\tilde{\varphi}_{ij}$ to $g_j\circ\tilde{\varphi}_{ij}$ for some $g_j\in G_j$ and satisfies (iv) (and automatically, (i) $\sim$ (iii)), then $\mathcal{S}$ and $\mathcal{S}'$ give the same orbifold structure. Two atlases give the same orbifold structure if their union is again a compatible atlas. We call $X$ the {\it underlying space} of the orbifold $M$, and denote it by the symbol $|M|$. Throughout this paper we assume $|M|$ to be paracompact.
\par
For each $x\in|M|$ a local chart $\varphi\circ r:\tilde{U}\rightarrow U$ of $M$ such that $U$ contains $x$ is called a {\it local chart around\/} $x$, and $\tilde{U}$ is called a {\it local uniformization around\/} $x$. The {\it local group at x}, denoted by $G_x$, is the isotropy group of any point in $\tilde{U}$ corresponding to $x$. This is well defined up to isomorphism. We call the order of $G_x$ the {\it index\/} (or {\it weight\/}) of $x$, and denote it by $w(x)$. The set $\{ x\in|M|\; |\; G_x\ne {id}\}$ is called the {\it singular set} of $M$ and denoted by $\Sigma M$. If $\Sigma M=\phi$, $M$ is a topological manifold.

A local chart $\tilde{U}_i\rightarrow \tilde{U}_i/G_i\rightarrow U$ around $x$ is {\it reduced\/} if $G_i$ is the local group at $x$. A point of $\Sigma M$ is called a {\it singular point}, and a point of $|M|-\Sigma M$ is called a {\it regular point}. A {\it base point\/} of $M$ is a fixed chosen point of $|M|-\Sigma M$. A {\it stratum\/} of $M$ is a maximal connected component of $|M|$ on which the orders of the local groups are constant. We denote the set of all $k$-dimensional strata by $\Sigma^{(k)}M$. An orbifold $M$ is {\it connected\/} if $|M|$ is connected.

\begin{remark}
If each $\tilde{\varphi}_{ij}$ is a smooth embedding in the above, we call $M$ a {\it smooth orbifold}. Note that (iv) is not assumed in the definition of orbifold (V-manifold) in \cite{Sa} or \cite{Th}.
\end{remark}

\section{Continuous map, homotopy and fundamental group}

In this section  first we review the definitions of (b-) continuous map between orbifolds. Isomorphisms and covering maps are regarded as special cases of continuous maps (cf. \cite{Th}).

\begin{definition}
Let $M$ and $N$ be orbifolds. By a {\it continuous map} $f:M\rightarrow N$, we shall mean a continuous map $|f|:|M|\rightarrow |N|$ together with a system $\mathcal{F}=(\{\varphi_i\circ r_i:\tilde{U}_i\rightarrow U_i\}_{i\in I},\{\psi_j\circ s_j:\tilde{V}_j\rightarrow V_j\}_{j\in J},\{\tilde{f}_{i\nu}:\tilde{U}_i\rightarrow \tilde{V}_{\nu}\}_{i\in I,\nu\in J_i})$ of a system of local charts of $M$, and of $N$, and a family of continuous maps, where $J_i=\{ \nu \in J\; |\; |f|(U_i)\subset V_{\nu }\}$, which satisfies the following:
\begin{enumerate}
\item[(0)] For any $i\in I$, $J_i\neq\emptyset$.
\item[(i)] For any $i\in I$ and any $\nu\in J_i$, $|f|\circ (\varphi_i\circ r_i)=(\psi_{\nu }\circ s_{\nu })\circ \tilde{f}_{i\nu }$.
\item[(ii)] For any $i\in I$, any $\nu\in J_i$, and any $\sigma_A\in$ Aut$(\tilde{U}_i,\varphi_i\circ r_i )$, there exists an element $\tau_A\in$ Aut$(\tilde{V}_{\nu },\psi_{\nu }\circ s_{\nu })$ such that $\tilde{f}_{i\nu }\circ\sigma_A =\tau_A\circ\tilde{f}_{i\nu }$.
\item[(iii)] If $U_i\subset U_j$, $|f|(U_i)\subset V_{\nu}$, $|f|(U_j)\subset V_{\mu}$, and $V_{\nu}\subset V_{\mu}$, then $\tilde{\psi}_{\nu\mu}\circ\tilde{f}_{i\nu}=\tilde{f}_{j\mu}\circ\tilde{\varphi}_{ij}$.
\end{enumerate}
Since $\tilde{U}_i-${\it Sing} $(G_i)$ is connected, $\#(\varphi_i\circ r_i)^{-1}(x)=\#$Aut$(\tilde{U}_i,\varphi_i\circ r_i)$, $x\in U_i-(U_i\cap\Sigma M)$. Then, Aut$(\tilde{U}_i,\varphi_i\circ r_i)=G_i$.

We call $|f|$ and $\{ \tilde{f}_{i\nu }\}$ the {\it underlying map} and the {\it structure maps} of the continuous map $f$, respectively. We often denote the above $f$ by $f=(|f|,\{\tilde{f}_{i\nu }\})$. We call $\mathcal{F}$ an {\it atlas\/} of $f$.
\par
If $f'$ is obtained from $f$ by changing each $\tilde{f}_{i\nu}$ to $g_{\nu}\circ\tilde{f}_{i\nu}$ for some $g_{\nu}\in G_{\nu}$ and satisfies (iii) (and automatically, (0) $\sim$ (ii)), then $f$ and $f'$ give the same continuous map structure. Two atlases give the same continuous map structure if their union is again a compatible atlas.

A continuous map $f:M\rightarrow N$ is {\it b-continuous\/} if there exists a point $x\in|M|-\Sigma M$ such that $|f|(x)\in|N|-\Sigma N$. It was called an orbi-map in \cite{japanfinite,finite,least}, etc. A b-continuous map induces a homomorphism between the fundamental groups and local fundamental groups of orbifolds, see Lemma \ref{pi1-induced-homo-lem} and \cite{finite}.
\end{definition}

\begin{remark}
Let $M$, $N$ be smooth orbifolds. Then a continuous map $f:M\rightarrow N$ is called {\it smooth\/} if each $\tilde{f}_{i\nu}$ is smooth.
\end{remark}

\begin{definition}
Let $M$ and $N$ be orbifolds. A continuous map $f=(|f|,\{\tilde{f}_{i\nu}\}):M\rightarrow N$ is an {\it isomorphism\/} if it satisfies the following:
\begin{enumerate}
\item[(i)] $|f|:|M|\rightarrow|N|$ is a homeomorphism.
\item[(ii)] For each $x\in|M|$, there exist reduced local charts $\tilde{U}_x\rightarrow\tilde{U}_x/G_x\cong U_x$ and $\tilde{U}'_{|f|(x)}\rightarrow\tilde{U}'_{|f|(x)}/G'_{|f|(x)}\cong U'_{|f|(x)}$ around $x$ and $|f|(x)$, respectively, and an isomorphism $\eta_x:G_x\rightarrow G'_{|f|(x)}$ between the local groups of $x$ and $|f|(x)$, respectively, such that each structure map $\tilde{f}_x:\tilde{U}_x\rightarrow\tilde{U}'_{|f|(x)}$ is a homeomorphism and for any $\sigma\in G_x$ and any $z\in\tilde{U}_x$, $\tilde{f}_x(\sigma z)=\eta_x(\sigma )\tilde{f}_x(z)$.
\end{enumerate}
An isomorphism is automatically a b-continuous map. If $\Sigma M=\Sigma N=\emptyset$, then an isomorphism $f:M\rightarrow N$ is a usual homeomorphism. If we discuss about an isomorphism between smooth orbifolds, we require that each $\tilde{f}_x$ is a diffeomorphism.
\end{definition}

\begin{definition}
Let $M$ and $N$ be orbifolds. Let $f,g:M\rightarrow N$ be continuous maps. By a {\it homotopy\/} from $f$ to $g$ (relative to $A\subset M$) we shall mean a continuous map $H:M\times[0,1]\rightarrow N$ which satisfies the following:
\begin{enumerate}
\item[(i)] $H(\cdot,0)=f(\cdot)$.
\item[(ii)] $H(\cdot,1)=g(\cdot)$.
\item[$\bigl($ (iii)] $H(a,s)=f(a)=g(a) \quad\mbox{for any}\; a\in A\;\mbox{and any}\; s\in[0,1].\;\bigr)$
\end{enumerate}

If the above $f$ and $g$ are b-continuous, $H$ is called a b-{\it homotopy}. A b-homotopy was called an orbi-homotopy in \cite{japanfinite,finite,least}, etc.
\end{definition}

\begin{definition}
Let $M$ and $N$ be orbifolds. A continuous map $f:M\rightarrow N$ is a {\it homotopy equivalent map\/} from $M$ to $N$ if there exists a continuous map $g:N\rightarrow M$ such that $g \circ f$ is homotopic to $id_M$ and $f\circ g$ is homotopic to $id_N$.

If the above $f$, $g$, $f\circ g$ and $g\circ f$ are all b-continuous, $f$ is called a b-{\it homotopy equivalent map}.
\end{definition}

\begin{example}
A cyclical ballic 3-orbifold $B^3(n)$ with index $n$ is b-homotopy equivalent to a discal 2-orbifold $D^2(n)$ with the same index $n$.
\end{example}

\begin{example}
Let $f$ be a trivial continuous map from a discal 2-orbifold $D^2(n)$ to a point $x$, that is, a connected 0-dimensional manifold. Let $g$ be a continuous map from $x$ to $D^2(n)$ defined by $g(x)=c$, where $c$ is the singular point of $D^2(n)$. Then $g\circ f$ (resp. $f\circ g$) is homotopic to $id_{D^2(n)}$ (resp. $id_x$), and $D^2(n)$ is homotopy equivalent to $x$. Note that $f$ is b-continuous but $g$ is not. Thus the pair of $f$ and $g$ does not give the b-homotopy equivalence of $D^2(n)$ and $x$. Indeed, $D^2(n)$ is not b-homotopy equivalnet to $x$. This fact will be derived from Theorem \ref{b-homotopy-inv-th} and the result of Subsection 9.2, or from computing their fundamental groups $\pi_1(D^2(n))\cong\mathbb{Z}_n$ and $\pi_1(\{ x\})=1$. See \ref{pi1-def} for the definition of $\pi_1$.
\end{example}

\begin{definition}
Let $M$ and $N$ be orbifolds. A continuous map $p=(|p|,\{\tilde{p}_{i\nu}\}):M\rightarrow N$ is a {\it covering map} if it satisfies the following:
\begin{enumerate}
\item[(i)] The underlying map $|p|:|M|\rightarrow|N|$ of $p$ is onto.
\item[(ii)] For each $x\in|N|$ and each $\tilde{x}\in|p|^{-1}(x)$, there exist reduced local charts $\tilde{U}_x\rightarrow\tilde{U}_x/G_x\cong U_x$ and $\tilde{U}_x\rightarrow\tilde{U}_x/G_{x,j}\cong V_{x,j}$, around $x$ and $\tilde{x}$, respectively, such that each structure map $\tilde{p}_x:\tilde{U}_x\rightarrow \tilde{U}_x$ is a homeomorphism and the following diagram commutes where $G_{x,j}$ is a subgroup of $G_x$, $V_{x,j}$ is the component of $|p|^{-1}(U_x)$ containing $\tilde{x}$ and $q$ is the natural projection:
\begin{equation}
\CD
\tilde{U}_x @>>> \tilde{U}_x/G_{x,j} @>>> V_{x,j} \\
@V\tilde{p}_xVV  @VqVV @V |p| VV \\
\tilde{U}_x @>>> \tilde{U}_x/G_x @>>> U_x
\endCD
\end{equation}
\end{enumerate}
\end{definition}

An isomorphism is a covering. An orbifold $M$ is {\it good\/} if there exists a covering $p:\tilde{M}\rightarrow M$ such that $\tilde{M}$ is a manifold.

A suborbifold $V$ of an orbifold $M$ is called the {\it regular neighborhood} of a point $x$ in $M$, if the underlying space $|V|$ is included in $U$ where $\varphi\circ r:\tilde{U}\rightarrow U$ is a reduced local chart around $x$ and $(\varphi\circ r)^{-1}(|V|)$ is the regular neighborhood of $(\varphi\circ r)^{-1}(x)$. A suborbifold $N$ is called the {\it regular neighborhood} of a subspace $Q$ in $M$, if $|N|$ is the regular neighborhood of $|Q|$ in $|M|$ and $N=\cup_{x\in |Q|}V_x$, where $V_x$ is the regular neighborhood of $x$ in $M$.
\begin{definition}\label{loop-def}
Let $M$ be an orbifold. A b-continuous map $a=(|a|,\{\tilde{a}_{i\nu}\}):[0,1]\rightarrow M$ with $|a|(0)\in|M|-\Sigma M$ is called a {\it path\/} in $M$. The point $|a|(0)$ is called the {\it initial point\/} of $a$, and $|a|(1)$ is called the {\it terminal point\/} of $a$. If $|a|(0)=p$ and $|a|(1)=q$, $a$ is called a path {\it from\/} $p$ {\it to\/} $q$. The initial and terminal points of $a$ are called the {\it end points\/} of $a$. If a path $a$ in $M$ satisfies that $|a|(0)=|a|(1)$, it is called a {\it loop\/} in $M$, and $|a|(0)\bigl(=|a|(1)\bigr)$ is called the {\it base point\/} of $a$. The product $a\cdot b$ of paths $a$ and $b$ can be defined if and only if $|a|(1)=|b|(0)$ (then $|a|(1)\in|M|-\Sigma M$), and the definition is as follows:
\begin{equation}
a\cdot b(t)=
\left\{
\begin{array}{ll}
a(2t)   & 0\leq t\leq\frac{1}{\; 2\;}, \\
\noalign{\vskip2mm}
b(2t-1) & \frac{1}{\; 2\;}\leq t\leq 1.
\end{array}\right.
\end{equation}
An orbifold $M$ is {\it arcwise connected\/} if for each two points $x,y\in|M|-\Sigma M$ there exists a path from $x$ to $y$ in $M$.
\end{definition}

\begin{definition}\label{pi1-def}
Let $M$ be an orbifold and let $x_0$ be a point of $|M|-\Sigma M$. The fundamental group $\pi_1(M,x_0)$ of $M$ with base point $x_0$ is the group of homotopy classes, relative to $\partial[0,1]$, of loops in $M$ with base point $x_0$, where the product of $[a]$ and $[b]$ is defined by
\begin{equation}
[a][b]:=[a\cdot b].
\end{equation}
\end{definition}

\begin{proposition}
Let $M$ be an arcwise connected orbifold and let $x$, $y$ be any two points of $|M|-\Sigma M$. Then the fundamental groups $\pi_1(M,x)$ and $\pi_1(M,y)$ are isomorphic.
\end{proposition}

\begin{proof}
Since $M$ is arcwise connected, there exists a path $\gamma:[0,1]\rightarrow M$ from $x$ to $y$. By using $\gamma$ we can construct an isomorphism from $\pi_1(M,x)$ to $\pi_1(M,y)$ as usual way.
\end{proof}

\begin{proposition}
Let $M$ be an orbifold with base point $x_0\in|M|-\Sigma M$, and let $p:\tilde{M}\rightarrow M$ be the universal covering of $M$. Note that $\tilde{M}$ is not necessarily a manifold. Then the fundamental group $\pi_1(M,x_0)$ of $M$ is isomorphic to the deck transformation group {\rm Aut}$(\tilde{M},p)$ of $p$.
\end{proposition}

Note that Thurston defined the fundamental group of an orbifold $M$ by Aut$(\tilde{M},p)$ in \cite{Th}.

\begin{lemma}\label{pi1-induced-homo-lem}
Let $M$, $N$ be orbifolds, and $f:M\rightarrow N$ be b-continuous map.
\begin{enumerate}
\item[(i)] The homomorphism between the fundamental groups of $M$ and $N$ can be induced by $f$.
\item[(ii)] Let $\varphi\circ r:\tilde{U}\rightarrow U$, $\psi\circ s:\tilde{V}\rightarrow V$ be local charts of $M$ and $N$, respectively, such that $|f|(U)\subset V$. Then the homomorphism between the local fundamental groups of $M|U$ and $N|V$ can be induced by $f$.
\end{enumerate}
\end{lemma}

\begin{proof}
(i) \quad Since $f$ is b-continuous, there exists a point $x_0\in |M|-\Sigma M$ such that $|f|(x_0)\in |N|-\Sigma N$. The induced homomorphism $f_*:\pi_1(M,x_0)\rightarrow \pi_1(N,|f|(x_0))$ can be defined by $f_*([a]):=[f\circ a]$ where $a$ is a loop with base point $x_0$. See \cite{finite} for the case of good orbifolds.

\noindent
(ii) \quad Let $y_0\in U-\Sigma (M|U)$ be a point and $c$ be a loop in $M|U$ with base point $y_0$. Take a path $\ell$ from $x_0$ to $y_0$ which does not intersect any singular point. Though $f\circ c$ might map into the singular set, $f\circ\ell$ intersects at least one regular point in $N|V$. That is, there exists $t_0\in[0,1]$ such that $|\ell|(t_0)$, $|f|\circ|\ell|(t_0)$ are regular points of $M|U$ and $N|V$, respectively. Thus we can define $f_*:\pi_1(M|U,|\ell|(t_0))\rightarrow \pi_1(N|V,|f|\circ|\ell|(t_0))$ similarly to the above.
\end{proof}

\section{Definition of t-singular homology}

Recall that a $q$-dimensional standard simplex is the subspace of $\mathbb{R}^{q+1}$ defined by (4.1). We denote it by $\triangle^q$.
\begin{equation}
\triangle^q=\Biggl\{x=(x_0,x_1,\dots,x_q)\in\mathbb{R}^{q+1}\;\Biggm|\;\sum_{i=0}^qx_i=1,\; x_i\geq 0,\; i=0,1,\dots,q\Biggr\}.
\end{equation}

\begin{definition}
Let $M$ be an orbifold.
\begin{enumerate}
\item[(i)] A $q$-dimensional {\it singular simplex\/} of $M$ is defined to be a continuous map $\Psi:\triangle^q\rightarrow M$.
\item[(ii)] A singular simplex $\Psi:\triangle^q\rightarrow M$ is {\it tame\/} if $|\Psi|^{-1}(\Sigma M)$ is finitely triangulable.
\end{enumerate}
\end{definition}

\begin{definition}
A $q$-dimensional singular simplex $\Psi:\triangle^q\rightarrow M$ of an orbifold $M$ is {\it non-transverse\/} with respect to an $i$-dimensional stratum $s\in\Sigma^{(i)}M$ if either $\Psi$ is not tame, or for a connected component $P$ of the inverse image $|\Psi|^{-1}(\bar{s})$ of the closure $\bar{s}$ of $s$ and an open neighbourhood $U\subset\triangle^q$ of $P$, there exists a homotopy $H:\triangle^q\times[0,1]\rightarrow M$ which satisfies the following:
\begin{enumerate}
\item[(i)] $H(x,0)=\Psi(x)$ \quad for any $x\in\triangle^q$.
\item[(ii)] $H(x,t)=\Psi(x)$ \quad for any $x\in\triangle^q-U$, any $t\in[0,1]$.
\item[(iii)] $H(x,t)\subset \cup_{z\in\Psi(P)}\stackrel{\circ}{N}_z$ \quad for any $x\in U$, any $t\in[0,1]$
\quad where $N_z$ is a regular neighbourhood of $z$ and $\stackrel{\circ}{N}_z$ is the interior of $N_z$.
\item[(iv)] $|H|(U,1)\cap\bar{s}=\emptyset$.
\end{enumerate}
Otherwise $\Psi$ is called {\it transverse\/} with respect to $s$. If for any $i$, $0\leq i\leq$ dim $\Sigma M$, and any $s\in\Sigma^{(i)}M$, $\Psi$ is transverse with respect to $s$, then $\Psi$ is called {\it transverse\/} with respect to $\Sigma M$.
\end{definition}

\begin{definition}\label{t-sing-def}
A $q$-dimensional {\it t-singular simplex\/} of an orbifold $M$ is a $q$-dimensional singular simplex $\Psi:\triangle^q\rightarrow M$ of $M$ such that each face of $\triangle^q$ is transverse with respect to $\Sigma M$.
\end{definition}

\begin{definition}
A $q$-dimensional {\it singular chain\/} (resp. {\it t-singular chain\/})  of an orbifold $M$ is a finite linear combination $\sum_jn_j\Psi_j$ of $q$-dimensional singular simplices (resp. t-singular simplices) of $M$. We denote the free abelian group with basis of all $q$-dimensional singular simplices (resp. t-singular simplices) of $M$ by $C_q(M)$ (resp. $t$-$C_q(M)$). By Definition \ref{t-sing-def}, $t$-$C_n(M)$ is a subgroup of $C_n(M)$.
\end{definition}

\begin{remark}
If an orbifold $M$ is a manifold, then the above $C_q(M)$ (resp. $t$-$C_q(M)$) coincides with the group of usual $q$-dimensional singular chains of $M$.
\end{remark}

The boundary operator is defined by the usual way as follows:

\begin{definition}
Let $M$ be an orbifold and $\Psi\in C_q(M)$. For $i=0,1,\dots,q$ we define $\partial_i\Psi\in C_{q-1}(M)$ by
\begin{equation}
\partial_i\Psi(x_0,x_1,\dots,x_{q-1})=\Psi(x_0,x_1,\dots,x_{i-1},0,x_i,\dots,x_{q-1}).
\end{equation}
Then the boundary operator $\partial:C_q(M)\rightarrow C_{q-1}(M)$ is defined by
\begin{equation}
\partial=\sum_{i=0}^q(-1)^i\partial_i,
\end{equation}
which satisfies that $\partial\circ\partial=0$. By Definition \ref{t-sing-def} it holds that $\partial(t$-$C_q(M))\subset t$-$C_{q-1}(M)$, and $t$-$C_*(M)$ is a subchain complex of $C_*(M)$.
\end{definition}

\begin{definition}
Let $M$ be an orbifold. As the usual way, we denote some subgroups of $C_q(M)$ as follows:
\begin{equation}
\begin{array}{ll}
Z_q(M)&=\{ c\in C_q(M) \;|\; \partial(c)=0 \}, \\
t\mbox{-}Z_q(M)&=\{c\in t\mbox{-}C_q(M)\;|\; \partial(c)=0 \}, \\
B_q(M)&=\{ c\in C_q(M) \;|\; \exists d\in C_{q+1}(M) \;\mbox{s.t.}\; \partial(d)=c \}, \\
t\mbox{-}B_q(M)&=\{ c\in t\mbox{-}C_q(M) \;|\; \exists d\in t\mbox{-}C_{q+1}(M) \;\mbox{s.t.}\; \partial(d)=c \}.
\end{array}
\end{equation}
Each element of $Z_q(M)$ (resp. $t$-$Z_q(M)$) is called a $q$-dimensional {\it cycle\/} (resp. {\it t-cycle}). And each element of $B_q(M)$ (resp. $t$-$B_q(M)$) is called a $q$-dimensional {\it boundary cycle\/} (resp. {\it t-boundary cycle\/}). Then the {\it singular homology\/} (resp. {\it t-singular homology\/}) of $M$ is defined by $Z_q(M)/B_q(M)$ (resp. $t$-$Z_q(M)/t$-$B_q(M)$) and is denoted by $H_q(M)$ (resp. $t$-$H_q(M)$).
\end{definition}

\begin{remark}
\begin{enumerate}
\item[(i)] Any continuous map between orbifolds induces a chain map between the singular complices. Thus the singular homology group is a homotopy invariant of orbifolds.
\item[(ii)] Any continuous map that is an embedding of orbifolds induces a chain map between the t-singular complices. In particular an isomorphism $f:M\rightarrow N$ between orbifolds induces a chain map $f_{\#}:t$-$C_*(M)\rightarrow t$-$C_*(N)$, which maps isomorphically $t$-$C_*(M)$ (resp. $t$-$Z_*(M)$, $t$-$B_*(M)$) to $t$-$C_*(N)$ (resp. $t$-$Z_*(N)$, $t$-$B_*(N)$). Thus the t-singular homology is an isomorphism invariant of orbifolds. In fact, the t-singular homology group is a b-homotopy invariant of orbifolds by Theorem \ref{b-homotopy-inv-th}. See Section 5 for the case of a b-continuous map between orbifolds that is not necessarily an embedding of orbifolds.
\end{enumerate}
\end{remark}

\section{t-modification}

Though any continuous (resp. b-continuous) map between orbifolds induces a homomorphism between the singular chain complices, it does not induce a homomorphism between the t-singular chain complices in general. This is caused by the fact that the composition of a t-singular simplex and a continuous (resp. b-continuous) map is not necessarily transverse with respect to the singular set of the orbifold. In this section we develope a \lq t-modification' of a singular chain and define a homomorphism between the t-singular homology groups induced by a b-continuous map between orbifolds.

\begin{definition}\label{sunny-def}
Let $E=\partial(\triangle^q\times[0,1])$- Int$(\triangle^q\times 1)$. Let $\ell$ be the line through the barycenters of $\triangle^q\times 0$ and $\triangle^q\times 1$. Let ${\bf O}$ be the point on $\ell$ such that $d({\bf O},\triangle^q\times 1)=1$ and ${\bf O}\not\in\triangle^q\times 0$. We define a map $\eta:\triangle^q\times[0,1]\rightarrow E$ be the intersection point of the line from ${\bf O}$ to $(x,s)\in\triangle^q\times[0,1]$ with $E$. That is, $\eta(x,s)={\bf O}z\cap E$, where $z=(x,s)$. Let $\phi:E\rightarrow M$ be a continuous map. Let $\{U_i\}$ be an open covering of $E$ which gives the orbifold structure of $E$ and ${\bf O}*U_i=\{t{\bf O}+(1-t)x\in\mathbb{R}^{q+1} \;|\; x\in U_i,0\leq t\leq 1 \}$. Let $\tilde{\phi}_{i\nu}:U_i\rightarrow\tilde{V}_j$ be the structure maps of $\phi$. Let $W_{ij}=({\bf O}*U_i)\cap\{(U_j\cap\triangle^q)\times[0,1]\}$. We define the {\it sunny extension\/} $\Phi:\triangle^q\times[0,1]\rightarrow M$ of $\phi$ as follows:
$|\Phi|:\triangle^q\times[0,1]\rightarrow M$ is defined as $|\Phi|(x,s)=|\phi|(\eta(x,s))$. $\tilde{\Phi}_{ij,\nu}:W_{ij}\rightarrow\tilde{V}_{\nu}$ is defined as $\tilde{\Phi}_{ij,\nu}(x,s)=\tilde{\phi}_{i\nu}(\eta(x,s))$.
\end{definition}

By $\lambda^{b}_{a_1\cdots a_r}$ we shall mean the inverse mapping of the embedding from the local uniformization $\tilde{V}_{a_1\cdots a_r}$  of $V_{a_1\cdots a_r}:=\cap_iV_{a_i}$ into $\tilde{V}_b$, $b\in\{a_1,\dots,a_r\}$.

\begin{lemma}\label{sunny1-lem}
Let $M$ be an orbifold and $\psi:\triangle^q=[x_0\cdots x_q]\rightarrow M$ a $q$-dimensional singular simplex. Let $\phi^s_k:[x_0\cdots\check{x}_k\cdots x_q]\rightarrow M$, $k=0,1,\dots,q$, be homotopies with $\phi^0_k=\psi|[x_0\cdots\check{x}_k\cdots x_q]$. Let $\{\tilde{\psi}_{i\nu}:U_i\rightarrow \tilde{V}_{\nu}\}_{i\in I,\nu\in J_i}$ and $\{(\tilde{\phi}^s_k)_{i_{(k,s)}\nu_{(k,s)}}:U_{i_{(k,s)}}\rightarrow \tilde{V}_{\nu_{(k,s)}}\}_{i_{(k,s)}\in I_{(k,s)},\nu_{(k,s)}\in J_{i_{(k,s)}}}$ be structure maps of $\psi$ and $\phi^s_k$, respectively, such that
\begin{equation*}
\lambda^{\nu_{(k,0)}}_{\nu_{(k,0)}\nu}(\tilde{\phi}^0_k)_{i_{(k,0)}\nu_{(k,0)}}|(U_{i_{(k,0)}}\cap U_i)=\lambda^{\nu}_{\nu_{(k,0)}\nu}\tilde{\psi}_{i\nu}|(U_{i_{(k,0)}}\cap U_i),
\end{equation*}
$k=0,\dots,q$. Then there exists a homotopy $\psi^s:[x_0\cdots x_q]\rightarrow M$ such that $\psi^s|[x_0\cdots \check{x}_k\cdots x_q]=\phi^s_k$ if and only if
\begin{equation*}
(|\phi^s_k|)|[x_0\cdots \check{x}_k\cdots \check{x}_{\ell}\cdots x_q]=(|\phi^s_{\ell}|)|[x_0\cdots\check{x}_k\cdots\check{x}_{\ell}\cdots x_q]
\end{equation*}
and
\begin{equation*}
\begin{array}{l}
\lambda^{\nu_{(k,s)}}_{\nu_{(k,s)}\nu_{(\ell,s)}}\circ(\tilde{\phi}^s_k)_{i_{(k,s)}\nu_{(k,s)}}|(U_{i_{(k,s)}}\cap U_{i_{(\ell,s)}}) \\
\noalign{\vskip2mm}
\quad\quad\quad\quad\quad =\lambda^{\nu_{(\ell,s)}}_{\nu_{(k,s)}\nu_{(\ell,s)}}\circ(\tilde{\phi}^s_{\ell})_{i_{(\ell,s)}\nu_{(\ell,s)}}|(U_{i_{(k,s)}}\cap U_{i_{(\ell,s)}}).
\end{array}
\end{equation*}
\end{lemma}

\begin{proof}
One direction is clear. If $(\tilde{\phi}^s_k)_{i_{(k,s)}\nu_{(k,s)}}$'s satisfy the hypothesis, then $\psi\cup\{\phi^s_i\}_{i\in I}$ defines a continuous map from $\partial(\triangle^q\times[0,1])- $Int$(\triangle^q\times 1)$ to $M$. Hence, we can construct $\tilde{\psi}^s$ by using of the sunny extension.
\end{proof}

\begin{lemma}\label{sunny2-lem}
Let $M$ be an orbifold and $\psi:\triangle^q=[x_0\cdots x_q]\rightarrow M$ a $q$-dimensional singular simplex. Let $\phi^s_k:[x_0\cdots\check{x}_k\cdots x_q]\rightarrow M$, $k=0,1,\dots,q$, be homotopies with $\phi^0_k=\psi|[x_0\cdots\check{x}_k\cdots x_q]$. Then there exists a homotopy $\psi^s:[x_0\cdots x_q]\rightarrow M$ such that $\psi^s|[x_0\cdots \check{x}_k\cdots x_q]=\phi^s_k$ if and only if $\phi^s_k|[x_0\cdots \check{x}_k\cdots\check{x}_{\ell}\cdots x_q]=\phi^s_{\ell}|[x_0\cdots \check{x}_k\cdots\check{x}_{\ell}\cdots x_q]$ and there exist points $x^s_{k\ell}\in[x_0\cdots \check{x}_k\cdots\check{x}_{\ell}\cdots x_q]$ such that $|\phi^s_k|(x^s_{k\ell})(=|\phi^s_{\ell}|(x^s_{k\ell}))\in|M|-\Sigma M$ for $\forall s,\forall k,\forall \ell=0,1,2,\dots,q$, $k<\ell$.
\end{lemma}

\begin{proof}
One direction is clear. Since $\triangle^q$ and $[0,1]$ are compact, there exists a finite open covering $P_0,\dots,P_N$ of $[0,1]$ such that for each $n\in\{0,1,\dots,N\}$, there exists a system of local charts $\{U_{i_{(k,n)}}\}_{i_{(k,s)}\in I_{(k,n)}}$ of $[x_0\cdots \check{x}_k\cdots x_q]$ such that $|\phi^s_k|(U_{i_{(k,n)}})$, $s\in P_n$, is included in a local chart of $M$. Let $\{V_{\nu_{(k,n)}}\}_{\nu_{(k,n)}\in J_{i_{(k,n)}}}$ be local charts of $M$ such that $|\phi^s_k|(U_{i_{(k,s)}})\subset V_{\nu_{(k,n)}}$. Let $\{\tilde{\psi}_{i\nu}:U_i\rightarrow \tilde{V}_{\nu}\}_{i\in I,\nu\in J_i}$ and $\{(\tilde{\phi}^s_k)_{i_{(k,s)}\nu_{(k,s)}}:U_{i_{(k,n)}}\rightarrow \tilde{V}_{\nu_{(k,n)}}\}_{i_{(k,n)}\in I_{(k,n)},\nu_{(k,n)}\in J_{i_{(k,n)}}}$ be structure maps of $\psi$ and $\phi^s_k$, respectively, such that
\begin{equation*}
\lambda^{\nu_{(k,0)}}_{\nu_{(k,0)}\nu}(\tilde{\phi}^0_k)_{i_{(k,0)}\nu_{(k,0)}}|(U_{i_{(k,0)}}\cap U_i)=\lambda^{\nu}_{\nu_{(k,0)}\nu}\tilde{\psi}_{i\nu}|(U_{i_{(k,0)}}\cap U_i),
\end{equation*}
$k=0,\dots,q$. From the hypothesis, there exists an element $g_{\nu_{(k,\ell,0)}}$ of the local group of $V_{\nu_{(k,0)}}\cap V_{\nu_{(\ell,0)}}$, such that
\begin{equation*}
\begin{array}{l}
\lambda^{\nu_{(k,0)}}_{\nu_{(k,0)}\nu_{(\ell,0)}}
(\tilde{\phi}^s_k)_{i_{(k,s)}\nu_{(k,s)}}
|
(U_{i_{(k,0)}}\cap U_{i_{(\ell,0)}}) \\
\noalign{\vskip2mm}
\quad\quad\quad\quad\quad= g_{\nu_{(k,\ell,0)}}
\lambda^{\nu_{(\ell,0)}}_{\nu_{(k,0)}\nu_{(\ell,0)}}
(\tilde{\phi}^s_{\ell})_{i_{(\ell,s)}\nu_{(\ell,s)}}
|
(U_{i_{(k,0)}}\cap U_{i_{(\ell,0)}}),
\end{array}
\end{equation*}
for $\forall s\in P_0$, $\forall k<\forall \ell$.
If $x^0_{k\ell}\in U_{i_{(k,0)}}\cap U_{i_{(\ell,0)}}$, then
\begin{equation*}
\lambda^{\nu_{(k,0)}}_{\nu_{(k,0)}\nu_{(\ell,0)}}
(\tilde{\phi}^s_k)_{i_{(k,s)}\nu_{(k,s)}}
(x^0_{k\ell})
=
g_{\nu_{(k,\ell,0)}}
\lambda^{\nu_{(\ell,0)}}_{\nu_{(k,0)}\nu_{(\ell,0)}}
(\tilde{\phi}^s_{\ell})_{i_{(\ell,s)}\nu_{(\ell,s)}}
(x^0_{k\ell}),
\end{equation*}
$s\in P_0$, namely
\begin{equation*}
\lambda^{\nu_{(k,0)}}_{\nu_{(k,0)}\nu_{(\ell,0)}}
(\tilde{\phi}^0_k)_{i_{(k,0)}\nu_{(k,0)}}
(x^0_{k\ell})
=
g_{\nu_{(k,\ell,0)}}
\lambda^{\nu_{(\ell,0)}}_{\nu_{(k,0)}\nu_{(\ell,0)}}
(\tilde{\phi}^0_k)_{i_{(\ell,0)}\nu_{(\ell,0)}}
(x^0_{k\ell}).
\end{equation*}
Let $U_i$ be a local chart of $\triangle^q$ such that $U_i\cap U_{i_{(k,0)}}\cap U_{i_{(\ell,0)}}\ni x^0_{k\ell}$. Then
\begin{equation*}
\begin{array}{lll}
\lambda^{\nu_{(k,0)}}_{\nu_{(k,0)}\nu_{(\ell,0)}\nu}(\tilde{\phi}^0_k)_{i_{(k,0)}\nu_{(k,0)}}(x^0_{k\ell})
&=& \lambda^{\nu}_{\nu_{(k,0)}\nu_{(\ell,0)}\nu}\tilde{\psi}_{i\nu}(x^0_{k\ell}) \\
\noalign{\vskip2mm}
&=& \lambda^{\nu_{(\ell,0)}}_{\nu_{(k,0)}\nu_{(\ell,0)}\nu}(\tilde{\phi}^0_k)_{i_{(\ell,0)}\nu_{(\ell,0)}}(x^0_{k\ell}).
\end{array}
\end{equation*}
Hence
\begin{equation*}
\begin{array}{l}
\lambda^{\nu_{(k,0)}}_{\nu_{(k,0)}\nu_{(\ell,0)}}(\lambda^{\nu_{(k,0)}}_{\nu_{(k,0)}\nu_{(\ell,0)}\nu})^{-1}\lambda^{\nu}_{\nu_{(k,0)}\nu_{(\ell,0)}\nu}\tilde{\psi}_{i\nu}(x^0_{k\ell}) \\
\noalign{\vskip2mm}
\quad\quad\quad\quad\quad =g_{\nu_{(k,\ell,0)}}\lambda^{\nu_{(\ell,0)}}_{\nu_{(k,0)}\nu_{(\ell,0)}}(\lambda^{\nu_{(\ell,0)}}_{\nu_{(k,0)}\nu_{(\ell,0)}\nu})^{-1}\lambda^{\nu}_{\nu_{(k,0)}\nu_{(\ell,0)}\nu}\tilde{\psi}_{i\nu}(x^0_{k\ell}).
\end{array}
\end{equation*}
From (iv) of the definition of orbifold,
\begin{equation*}
\begin{array}{l}
\lambda^{\nu_{(k,0)}}_{\nu_{(k,0)}\nu_{(\ell,0)}}(\lambda^{\nu_{(k,0)}}_{\nu_{(k,0)}\nu_{(\ell,0)}\nu})^{-1}\lambda^{\nu}_{\nu_{(k,0)}\nu_{(\ell,0)}\nu}\tilde{\psi}_{i\nu}(x^0_{k\ell}) \\
\noalign{\vskip2mm}
\quad\quad\quad\quad\quad\quad\quad\quad =\lambda^{\nu_{(\ell,0)}}_{\nu_{(k,0)}\nu_{(\ell,0)}}(\lambda^{\nu_{(\ell,0)}}_{\nu_{(k,0)}\nu_{(\ell,0)}\nu})^{-1}\lambda^{\nu}_{\nu_{(k,0)}\nu_{(\ell,0)}\nu}\tilde{\psi}_{i\nu}(x^0_{k\ell}).
\end{array}
\end{equation*}
Furthermore, since $|\psi|(x^0_{k\ell})\in|M|-\Sigma M$,
\begin{equation*}
\lambda^{\nu_{(k,0)}}_{\nu_{(k,0)}\nu_{(\ell,0)}}(\lambda^{\nu_{(k,0)}}_{\nu_{(k,0)}\nu_{(\ell,0)}\nu})^{-1}\lambda^{\nu}_{\nu_{(k,0)}\nu_{(\ell,0)}\nu}\tilde{\psi}_{i\nu}(x^0_{k\ell})
\end{equation*}
is a non-singular point. Hence, $g_{\nu_{(k,\ell,0)}}=e$.

Suppose $x^0_{k\ell}\not\in U_{i_{(k,0)}}\cap U_{i_{(\ell,0)}}$. By using the path in $[x_0\cdots \check{x}_k\cdots\check{x}_{\ell}\cdots x_q]$ from $x^0_{k\ell}$ to $U_{i_{(k,0)}}\cap U_{i_{(\ell,0)}}$, we can construct local charts $U'_{i_{(k,0)}}$, $U'_{i_{(\ell,0)}}$ and structure maps $(\tilde{\phi}^s_k)'_{i_{(k,s)}\nu_{(k,s)}}:U'_{i_{(k,0)}}\rightarrow \tilde{V}_{\nu_{(k,0)}}$, $(\tilde{\phi}^s_{\ell})'_{i_{(\ell,s)}\nu_{(\ell,s)}}:U'_{i_{(\ell,0)}}\rightarrow \tilde{V}_{\nu_{(\ell,0)}}$ of $(\phi^s_k)'$, $(\phi^s_{\ell})'$, respectively, such that $U'_{i_{(k,0)}}\supset U_{i_{(k,0)}}$, $U'_{i_{(\ell,0)}}\supset U_{i_{(\ell,0)}}$, 
$(\tilde{\phi}^s_k)'_{i_{(k,s)}\nu_{(k,s)}}|U_{i_{(k,0)}}=(\tilde{\phi}^s_k)_{i_{(k,s)}\nu_{(k,s)}}$,
$(\tilde{\phi}^s_{\ell})'_{i_{(\ell,s)}\nu_{(\ell,s)}}|U_{i_{(\ell,0)}}=(\tilde{\phi}^s_{\ell})_{i_{(\ell,s)}\nu_{(\ell,s)}}$, and
$x^0_{k\ell}\in U'_{i_{(k,0)}}\cap U'_{i_{(\ell,0)}}$.
Then,we have
\begin{equation*}
\begin{array}{l}
(\lambda^{\nu_{(k,0)}}_{\nu_{(k,0)}\nu_{(\ell,0)}})'
(\tilde{\phi}^s_k)'_{i_{(k,s)}\nu_{(k,s)}}
|
(U'_{i_{(k,0)}}\cap U'_{i_{(\ell,0)}}) \\
\noalign{\vskip2mm}
\quad\quad\quad\quad\quad=
(\lambda^{\nu_{(\ell,0)}}_{\nu_{(k,0)}\nu_{(\ell,0)}})'
(\tilde{\phi}^s_{\ell})'_{i_{(\ell,s)}\nu_{(\ell,s)}}
|
(U'_{i_{(k,0)}}\cap U'_{i_{(\ell,0)}}).
\end{array}
\end{equation*}
By restricting $(\tilde{\phi}^s_k)'_{i_{(k,s)}\nu_{(k,s)}}$ and $(\tilde{\phi}^s_{\ell})'_{i_{(\ell,s)}\nu_{(\ell,s)}}$ to $U_{i_{(k,0)}}\cap U_{i_{(\ell,0)}}$ and using (iv) of the definition of orbifold, we have
\begin{equation*}
\begin{array}{l}
\lambda^{\nu_{(k,0)}}_{\nu_{(k,0)}\nu_{(\ell,0)}}
(\tilde{\phi}^s_k)_{i_{(k,s)}\nu_{(k,s)}}
|
(U_{i_{(k,0)}}\cap U_{i_{(\ell,0)}}) \\
\noalign{\vskip2mm}
\quad\quad\quad\quad\quad=
\lambda^{\nu_{(\ell,0)}}_{\nu_{(k,0)}\nu_{(\ell,0)}}
(\tilde{\phi}^s_{\ell})_{i_{(\ell,s)}\nu_{(\ell,s)}}
|
(U_{i_{(k,0)}}\cap U_{i_{(\ell,0)}}),
\end{array}
\end{equation*}
Similarly,
\begin{equation*}
\begin{array}{l}
\lambda^{\nu_{(k,1)}}_{\nu_{(k,1)}\nu_{(\ell,1)}}
(\tilde{\phi}^s_k)_{i_{(k,s)}\nu_{(k,s)}}
|
(U_{i_{(k,1)}}\cap U_{i_{(\ell,1)}}) \\
\noalign{\vskip2mm}
\quad\quad\quad\quad\quad=
\lambda^{\nu_{(\ell,1)}}_{\nu_{(k,1)}\nu_{(\ell,1)}}
(\tilde{\phi}^s_{\ell})_{i_{(\ell,s)}\nu_{(\ell,s)}}
|
(U_{i_{(k,1)}}\cap U_{i_{(\ell,1)}}),
\end{array}
\end{equation*}
$s\in P_1$, is proved by the fact that
\begin{equation*}
\begin{array}{l}
\lambda^{\nu_{(k,1)}}_{\nu_{(k,1)}\nu_{(\ell,1)}}
(\tilde{\phi}^s_k)_{i_{(k,\varepsilon)}\nu_{(k,\varepsilon)}}
|
(x^{\varepsilon}_{k\ell}) \\
\noalign{\vskip2mm}
\quad\quad\quad\quad\quad=
\lambda^{\nu_{(\ell,1)}}_{\nu_{(k,1)}\nu_{(\ell,1)}}
(\tilde{\phi}^s_{\ell})_{i_{(\ell,\varepsilon)}\nu_{(\ell,\varepsilon)}}
|
(x^{\varepsilon}_{k\ell}),
\end{array}
\end{equation*}
$\varepsilon\in P_0\cap P_1$. Inductively, $s\in P_2,\dots,P_N$, we have
\begin{equation*}
\begin{array}{l}
\lambda^{\nu_{(k,m)}}_{\nu_{(k,m)}\nu_{(\ell,m)}}
(\tilde{\phi}^s_k)_{i_{(k,s)}\nu_{(k,s)}}
|
(U_{i_{(k,m)}}\cap U_{i_{(\ell,m)}}) \\
\noalign{\vskip2mm}
\quad\quad\quad\quad\quad=
\lambda^{\nu_{(\ell,m)}}_{\nu_{(k,m)}\nu_{(\ell,m)}}
(\tilde{\phi}^s_{\ell})_{i_{(\ell,s)}\nu_{(\ell,s)}}
|
(U_{i_{(k,m)}}\cap U_{i_{(\ell,m)}}),
\end{array}
\end{equation*}
$m=2,3,\dots,N$, for $\forall s$. Then, by Lemma \ref{sunny1-lem}, we have the conclusion.
\end{proof}

\begin{lemma}
For any singular simplex $\varphi:\triangle^q\rightarrow M$, there exists a homotopy $\varphi^s:\triangle\rightarrow M$ such that $\varphi^0=\varphi$, $\varphi^s|\triangle^{(0)}=\varphi|\triangle^{(0)}$, and $\varphi^1$ is tame. Furthermore, if $\varphi|\partial\triangle^q$ is transverse with respect to $\Sigma M$, then $\varphi^s|\partial\triangle^q=\varphi|\partial\triangle^q$.
\end{lemma}

\begin{proof}
We can construct $\varphi^s$ by using of a PL-approximation of $\varphi$.
\end{proof}

Let $\triangle^q_{(i_0\cdots i_k)}$ be the $k$-dimensional face of $\triangle^q$ defined by
\begin{equation}
\triangle^q_{(i_0\cdots i_k)}=\{(x_0,\dots,x_q)\in\triangle^q\:|\; x_{i_0}+\cdots+x_{i_k}=1\}
\end{equation}
and let $\varepsilon^q_{(i_0\cdots i_k)}:\triangle^k\rightarrow\triangle^q_{(i_0\cdots i_k)}$ be the linear homeomorphism define by $\varepsilon^q_{(i_0\cdots i_k)}(x_j):=x_{i_j}$.

\begin{definition}\label{t-mod-def}
By a {\it t-modification\/} of a $q$-dimensional singular chain $c=\sum^{\ell}_{j=1}n_j\Psi^j$ of an orbifold $M$ we shall mean a t-singular chain $c^T=\sum^{\ell}_{j=1}n_j(\Psi^j)^T$ such that there exist homotopies $\Psi^1_s,\dots,\Psi^{\ell}_s$ ($s\in[0,1]$) that satisfy the following:
\begin{enumerate}
\item[(i)] $\Psi^j_0=\Psi^j$, $\Psi^j_1=(\Psi^j)^T$ \quad for $\forall j=1,2,\dots,\ell$.
\item[(ii)] If
\begin{equation}
(\Psi^i|\triangle^q_{(i_0\cdots i_k)})\circ\varepsilon^q_{(i_0\cdots i_k)}=(\Psi^j|\triangle^q_{(j_0\cdots j_k)})\circ\varepsilon^q_{(j_0\cdots j_k)},
\end{equation}
for some $i$, $j$, then it holds that
\begin{equation}
(\Psi^i_s|\triangle^q_{(i_0\cdots i_k)})\circ\varepsilon^q_{(i_0\cdots i_k)}=(\Psi^j_s|\triangle^q_{(j_0\cdots j_k)})\circ\varepsilon^q_{(j_0\cdots j_k)},
\quad \mbox{for} \;\forall \; s\in[0,1].
\end{equation}
\end{enumerate}
We call $(c^T)_s=\sum n_j\Psi^j_s$ ($s\in[0,1]$) a {\it t-homotopy\/} of $c$.
\end{definition}

\begin{remark}
By (ii) of Definition \ref{t-mod-def} a t-modification of a singular chain has to preserve the face conditions. The following singular chain (simplex) $\Psi$ has no t-modifications, see Figure 5.1. Note that neither of the singular chains (simplices) $\Psi_1$ and $\Psi_2$ illustrated by Figures 5.2 and 5.3 is not a t-modification of $\Psi$.
\end{remark}

\begin{figure}[htbp]
\begin{center}
\includegraphics[width=.5\linewidth]{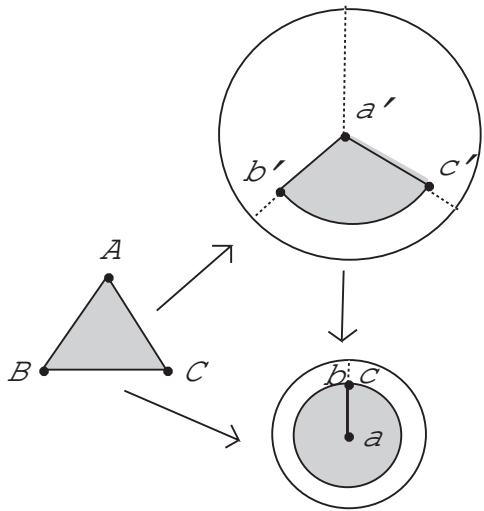}
\end{center}
\caption{$\Psi=(|\Psi|,\tilde{\Psi}):\triangle^2\rightarrow D^2(3)$, $\tilde{\Psi}(A)=a',\tilde{\Psi}(B)=b',\tilde{\Psi}(C)=c',|\Psi|(A)=a,|\Psi|(B)=b,|\Psi|(C)=c$}
\end{figure}

\begin{figure}[htbp]
\begin{center}
\includegraphics[width=.5\linewidth]{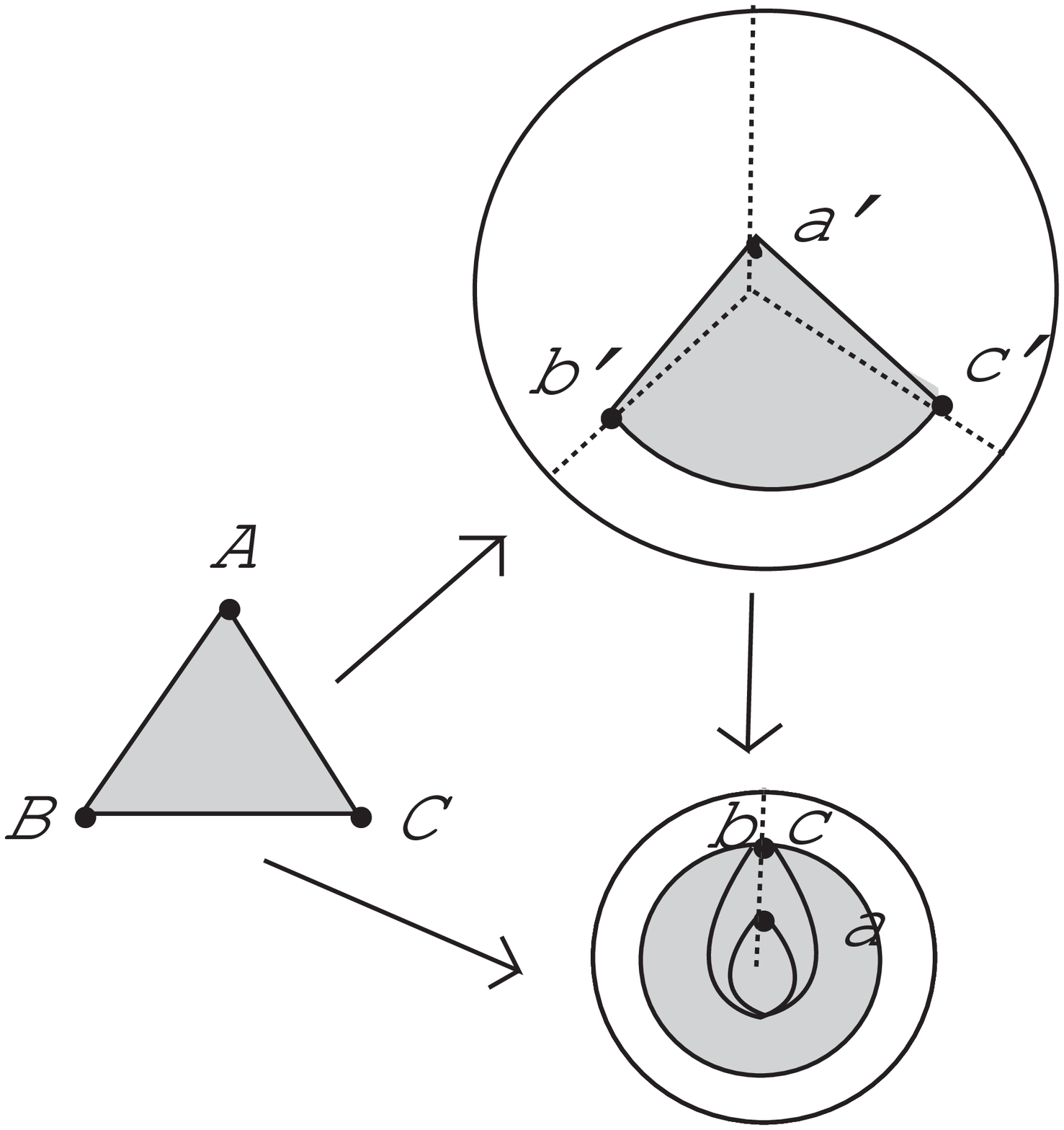}
\end{center}
\caption{$\Psi_1:\triangle^2\rightarrow D^2(3)$}
\end{figure}

\begin{figure}[htbp]
\begin{center}
\includegraphics[width=.5\linewidth]{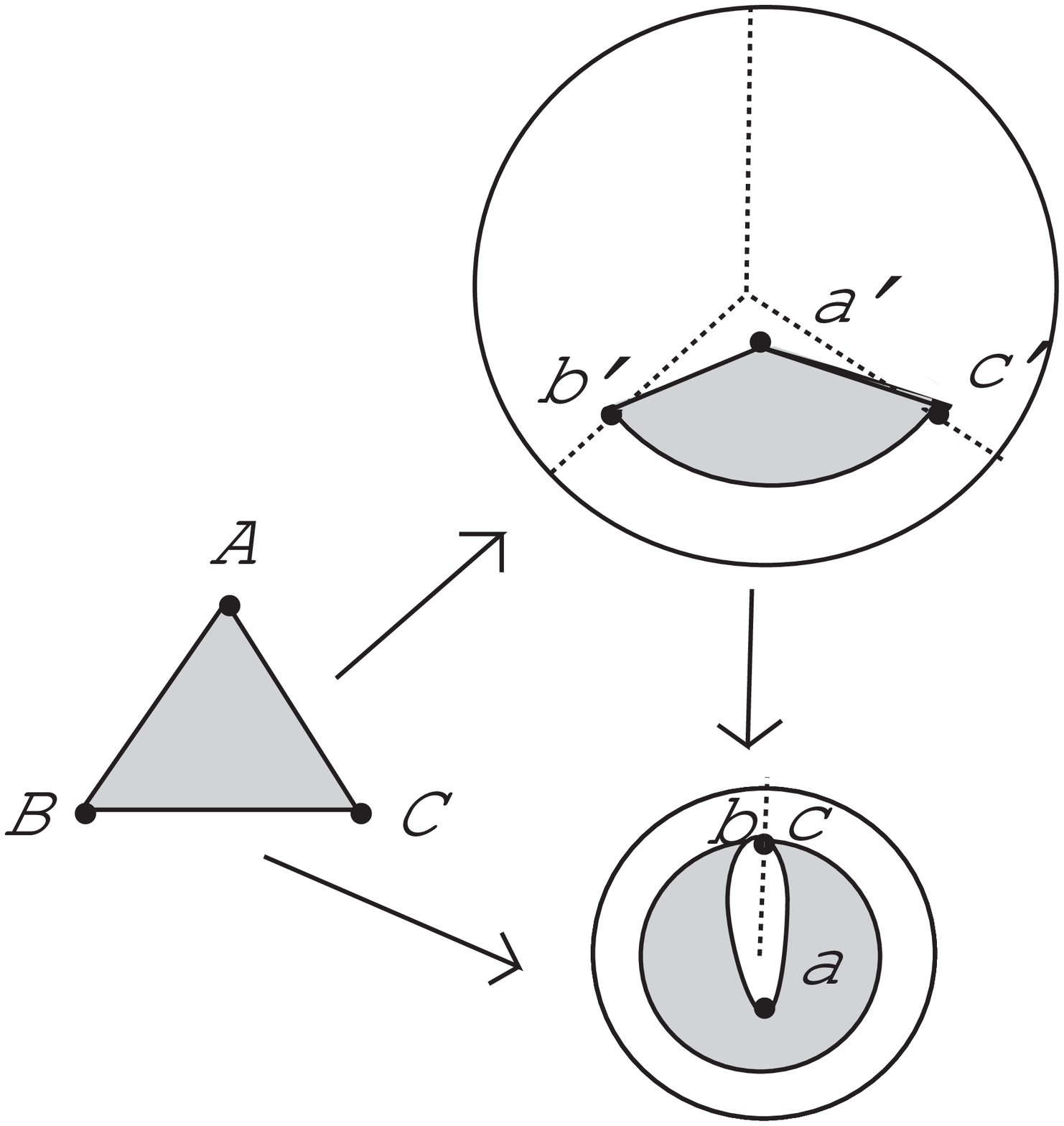}
\end{center}
\caption{$\Psi_2:\triangle^2\rightarrow D^2(3)$}
\end{figure}

\begin{definition}
A $q$-dimensional singular simplex $\Psi=(|\Psi|,\{\tilde{\Psi}_{ij}\})$ of $M$ is {\it pre-transverse\/} if the underlying map $|\Psi|$ of $\Psi$ maps all the vertices of $\triangle^q$ into $|M|-\Sigma M$. A $q$-dimensional singular chain $c=\sum m_i\Psi^i$ of $M$ is {\it pre-transverse\/} if each singular simplex $\Psi^i$ is pre-transverse.
\end{definition}

Note that every t-singular chain is pre-transverse.

\begin{lemma}\label{t-mod-lem}
Let $M$ be an orbifold and let $c=\sum^{\ell}_{\mu=1}n_{\mu}\varphi^{\mu}$ be a $q$-dimensional singular chain of $M$. Let $c'$ be a sum of faces of $c$ all of which are already transverse with respect to $\Sigma M$. If $c$ is pre-transverse, then there exists a t-homotopy $c_s$ of $c$ fixing $c'$.
\end{lemma}

\begin{proof}
To construct a t-homotopy $c_s=\sum^{\ell}_{\mu=1}n_{\mu}\varphi^{\mu}_s$ of $c$ we will define
\begin{equation}
(\varphi^{\mu}_{(i_0\cdots i_k)})_s:\triangle^q_{(i_0\cdots i_k)}\rightarrow M,\; \mu=1,2,\dots,\ell,\; 0\leq i_0<\cdots<i_k\leq q
\end{equation}
inductively on $k$ such that
\begin{enumerate}
\item[(1)] $(\varphi^{\mu}_{(i_0\cdots i_k)})_0=\varphi^{\mu}|\triangle^q_{(i_0\cdots i_k)}$.
\item[(2)] $(\varphi^{\mu}_{(i_0\cdots i_k)})_1$ is transverse with respect to $\Sigma M$.
\item[(3)] $(\varphi^{\mu}_{(i_0\cdots i_k)})_s$'s satisfy the face condition of $\varphi^{\mu}_{(i_0\cdots i_k)}$'s.
\item[(4)] $(\varphi^{\mu}_{(i_0\cdots i_k)})_s|\triangle^q_{(i_0\cdots \check{i}_a\cdots i_k)}=(\varphi^{\mu}_{(i_0\cdots \check{i}_a\cdots i_k)})_s$.
\item[(5)] If $\varphi^{\mu}|\triangle^q_{(i_0\cdots i_k)}$ is transverse with respect to $\Sigma M$, then $(\varphi^{\mu}_{(i_0\cdots i_k)})_s=\varphi^{\mu}|\triangle^q_{(i_0\cdots i_k)}$.
\end{enumerate}
If $k=0$, put $(\varphi^{\mu}_{(i_0)})_s:=\varphi^{\mu}|\triangle^q_{(i_0)}$. It is transverse with respect to $\Sigma M$, since $c$ is pre-transverse.

Suppose that we have done for the ($k-1$)-skelton. We order $\{(\mu,i_0,\dots,i_k)\}$ lexicographically and will construct $(\varphi^{\mu}_{(i_0\cdots i_k)})_s$'s, $\mu=1,2,\dots,\ell$, $0\leq i_0<\dots<i_k\leq q$, in that order as follows:

1) \quad If $\varphi^{\mu}|\triangle^q_{(i_0\cdots i_k)}$ is transverse with respect to $\Sigma M$, then we define as
\begin{equation*}
(\varphi^{\mu}_{(i_0\cdots i_k)})_s:=\varphi^{\mu}|\triangle^q_{(i_0\cdots i_k)}.
\end{equation*}

2) \quad If $\varphi^{\mu}|\triangle^q_{(i_0\cdots i_k)}$ is not transverse with respect to $\Sigma M$ and for some $(\nu,j_0,\dots,j_k)<(\mu,i_0,\dots,i_k)$, $(\varphi^{\nu}|\triangle^q_{(j_0\cdots j_k)})\circ\varepsilon^q_{(j_0\cdots j_k)}=(\varphi^{\mu}|\triangle^q_{(i_0\cdots i_k)})\circ\varepsilon^q_{(i_0\cdots i_k)}$, then we define as
\begin{equation*}
(\varphi^{\mu}_{(i_0\cdots i_k)})_s=(\varphi^{\nu}_{(j_0\cdots j_k)})_s\circ\varepsilon^q_{(j_0\cdots j_k)}\circ(\varepsilon^q_{(i_0\cdots i_k)})^{-1}.
\end{equation*}

3) \quad Suppose $\varphi^{\mu}|\triangle^q_{(i_0\cdots i_k)}$ is not transverse with repsect to $\Sigma M$ and there exists no $(\nu,j_0,\dots,j_k)<(\mu,i_0,\dots,i_k)$, such that $(\varphi^{\nu}|\triangle^q_{(j_0\cdots j_k)})\circ\varepsilon^q_{(j_0\cdots j_k)}=(\varphi^{\mu}|\triangle^q_{(i_0\cdots i_k)})\circ\varepsilon^q_{(i_0\cdots i_k)}$. By the inductive assumption $(\varphi^{\mu}_{(i_0\cdots\check{i}_a\cdots i_k)})_s$, $a=0,1,2,\dots,k$, are already defined. Then, by Lemma \ref{sunny2-lem} there exists a homotopy $(\phi^{\mu}_{(i_0\cdots i_k)})_s:\triangle^q_{(i_0\cdots i_k)}\rightarrow M$ such that $(\phi^{\mu}_{(i_0\cdots i_k)})_s|\triangle^q_{(i_0\cdots\check{i}_a\cdots i_k)}=(\varphi^{\mu}_{(i_0\cdots\check{i}_a\cdots i_k)})_s$ and $(\phi^{\mu}_{(i_0\cdots i_k)})_0=\varphi^{\mu}|\triangle^q_{(i_0\cdots i_k)}$. Since $(\phi^{\mu}_{(i_0\cdots i_k)})_s$ is derived from the sunny extension, $(\phi^{\mu}_{(i_0\cdots i_k)})_1$ is not transverse with respect to $\Sigma M$. Then for a dimension $k_0$ and for a component $s_{k_0}\in\Sigma^{(k_0)}M$, there exist a component $P_0$ of $(\phi^{\mu}_{(i_0\cdots i_k)})_1^{-1}(\bar{s}_{k_0})$, an open neighbourhood $U_0\subset \triangle^q_{(i_0\cdots i_k)}$ of $P_0$ and a homotopy $\Phi^{k_0}_s:\triangle^q_{(i_0\cdots i_k)}\rightarrow M$ which satisfy the following:
\begin{enumerate}
\item[(1)] $\Phi^{k_0}_s(x)=(\phi^{\mu}_{(i_0\cdots i_k)})_1(x)$ for $\forall x\in \triangle^q_{(i_0\cdots i_k)}-U_0,\;\forall\; s\in[0,1]$.
\item[(2)] $\Phi^{k_0}_0(x)=(\phi^{\mu}_{(i_0\cdots i_k)})_1(x)$ for $\forall x\in \triangle^q_{(i_0\cdots i_k)}$.
\item[(3)] $|\Phi^{k_0}_s|(x)\in\cup_{z\in(\phi^{\mu}_{(i_0\cdots i_k)})_1(P_0)}\stackrel{\circ}{N}_z$ for $\forall x\in U_0,\;\forall\; s\in[0,1]$.
\item[(4)] $|\Phi^{k_0}_1|(U_0)\cap\bar{s}_{k_0}=\emptyset$.
\end{enumerate}
We may assume that $k_0$ is the minimum one satisfying the above. By the minimality of $k_0$, $(\phi^{\mu}_{(i_0\cdots i_k)})_1(P_0)\subset s_{k_0}$, that is, $(\phi^{\mu}_{(i_0\cdots i_k)})_1(P_0)\cap(\bar{s}_{k_0}-s_{k_0})=\emptyset$.
Thus
\begin{equation}\label{Nz-eq}
N_z\cap(\Sigma^{(k_0)}M\cup\cdots\cup\Sigma^{(0)}M)\subset\bar{s}_{k_0} \quad \mbox{for} \;\forall z\in(\phi^{\mu}_{(i_0\cdots i_k)})_1(P_0).
\end{equation}
By \eqref{Nz-eq} and (4),
\begin{equation}
|\Phi^{k_0}_1|(U_0)\subset|M|-(\Sigma^{(k_0)}M\cup\cdots\cup\Sigma^{(0)}M).
\end{equation}

If $\Phi^{k_0}_1$ is transverse with respect to $\Sigma M$, we define $(\varphi^{\mu}_{(i_0\cdots i_k)})_s$ as the product of $(\phi^{\mu}_{(i_0\cdots i_k)})_s$ and $\Phi^{k_0}_s$. If not, for an index $k_1$ and for a component $s_{k_1}'\in\Sigma^{(k_1)}M$ there exist a component $P_1$ of $(\Phi^{k_0}_1)^{-1}(\bar{s}_{k_1}')$, an open neighbourhood $U_1\subset \triangle^q_{(i_0\cdots i_k)}$ of $P_1$, and a homotopy $\Phi^{k_1}_s:\triangle^q_{(i_0\cdots i_k)}\rightarrow M$ which satisfy the following:
\begin{enumerate}
\item[(1')] $\Phi^{k_1}_s(x)=\Phi^{k_0}_1(x)$ for $\forall x\in \triangle^q_{(i_0\cdots i_k)}-U_1$, $\forall\; s\in[0,1]$.
\item[(2')] $\Phi^{k_1}_0(x)=\Phi^{k_0}_1(x)$ for $\forall x\in \triangle^q_{(i_0\cdots i_k)}$.
\item[(3')] $|\Phi^{k_1}_s|(x)\subset \cup_{z\in\Phi^{k_0}_1(P_1)}\stackrel{\circ}{N}_z$ for $\forall x\in U_1$, $\forall\; s\in[0,1]$.
\item[(4')] $|\Phi^{k_1}_1|(U_1)\cap\bar{s}_{k_1}'=\emptyset$.
\end{enumerate}
We may assume that $k_1$ is the minimum one satisfying the above.

a) If \quad $k_1=k_0$ and $s_{k_1}'=s_{k_0}$, the argument can be reduced to the case b) or c) since the number of the components of $(\phi^{\mu}_{(i_0\cdots i_k)})_1^{-1}(\bar{s}_{k_0})$ is finite.

b) \quad If $k_1=k_0$ and $s_{k_1}'\neq s_{k_0}$, the argument can be reduced to the case c) since the number of the components of $(\phi^{\mu}_{(i_0\cdots i_k)})_1(\triangle^q_{(i_0\cdots i_k)})\cap\Sigma^{(k_0)}M$ is finite.

c) If $k_1\neq k_0$, then $k_1>k_0$. Since $|\Phi^{k_0}_1|(P_1)\subset s_{k_1}'$, it holds that
\begin{equation}\label{Nz2-eq}
N_z\cap(\Sigma^{(k_1)}M\cup\cdots\cup\Sigma^{(0)}M)\subset\bar{s}_{k_1}'\quad\mbox{for}\;\forall z\in \Phi^{k_0}_1(P_1).
\end{equation}
By \eqref{Nz2-eq} and (4'), it holds that
\begin{equation}
|\Phi^{k_1}_1|(U_1)\subset|M|-(\Sigma^{(k_1)}M\cup\cdots\cup\Sigma^{(0)}M).
\end{equation}
If we repeat this process in finte times, then, for some $q$, $\Phi^{k_q}_1$ is transverse with respect to $\Sigma M$. We define $(\varphi^{\mu}_{(i_0\cdots i_k)})_s$ as the product of $(\phi^{\mu}_{(i_0\cdots i_k)})_s,\Phi^{k_0}_s,\dots,\Phi^{k_q}_s$.
Thus we have the desired homotopy.
\end{proof}

\begin{lemma}\label{tpartial-lem}
Let $c$ be a singular chain of an orbifold $M$. Suppose that $\partial c$ is transverse with respect to $\Sigma M$.  If $c_s$ is a t-homotopy of $c$ fixing $\partial c$, then for any $s\in[0,1]$ it holds that $\partial c_s=\partial c$. In particular, if $c$ is a cycle, then so is $c_s$.
\end{lemma}

\begin{proof}
Put $c=\sum_{i=1}^{\ell}n_i\varphi^i$. Its boundary is as follows:
\begin{equation}\label{partial-eq}
\begin{align}
\partial c
&=\sum_{i=1}^{\ell}n_i\partial\varphi^i \\
&=\sum_{i=1}^{\ell}\sum_{j=0}^n(-1)^jn_i\varphi^i\circ\varepsilon_j\nonumber \\
&=\sum_{(a,b)\in D}N_{(a,b)}\varphi^a\circ\varepsilon_b \nonumber
\end{align}
\end{equation}
where each $N_{(a,b)}$ is an integer, and 
\begin{equation}
D=\{(a,b)\in \{1,2,\dots,\ell\}\times\{0,1,2,\dots,n\}\;|\; N_{(a,b)}\neq 0\}.
\end{equation}
Let $c_s=\sum_{i=1}^{\ell}n_i\varphi^i_s$ be a t-homotopy of $c$. Since $c_s$ fixes $\partial c$,
\begin{equation}\label{varphi-eq}
\varphi^a_s\circ\varepsilon_b=\varphi^a\circ\varepsilon_b
\end{equation}
for any $(a,b)\in D$. By (ii) of Definition \ref{t-mod-def}, \eqref{varphi-eq}, and \eqref{partial-eq}
\begin{equation}
\begin{align}
\partial c_s
&=\sum_{i=1}^{\ell}\sum_{j=0}^n(-1)^jn_i\varphi^i_s\circ\varepsilon_j \\
&=\sum_{(a,b)\in D}N_{(a,b)}\varphi^a_s\circ\varepsilon_b \nonumber \\
&=\sum_{(a,b)\in D}N_{(a,b)}\varphi^a\circ\varepsilon_b \nonumber \\
&=\partial c. \nonumber
\end{align}
\end{equation}
\end{proof}

Note that if $c_s$ does not satisfy the face condition, then its boundary can be different from the above, that is,
$$\partial c_s=\sum_{(a,b)\in D}N_{(a,b)}\varphi^a_s\circ\varepsilon_b+B_s$$
where $B_s\in C_{q-1}(M)$.

\begin{definition}
Let $M$ and $N$ be orbifolds, $f:M\rightarrow N$ a b-continuous map, and $c=\sum_{j=1}^{\ell}n_j\varphi^j$ a $q$-dimensional pre-transversal singular chain of $M$. A {\it pre-t-modification\/} of a singular chain $f\circ c=\sum^{\ell}_{j=1}n_j(f\circ\varphi^j)$ with respect to $c$ is a pre-transversal singular chain $d=\sum^{\ell}_{j=1}n_j\Psi^j$ of $N$ such that there exist homotopies $\Psi^1_s,\dots,\Psi^{\ell}_s$ ($s\in[0,1]$) which satisfy the following:
\begin{enumerate}
\item[(i)] $\Psi^j_0=f\circ\varphi^j$, $\Psi^j_1=\Psi^j$ \quad for $\forall j=1,2,\dots,\ell$.
\item[(ii)] If for $\varphi^i$ and $\varphi^j$ there exist embeddings $\varepsilon^q_{i_0\cdots i_{q-k-1}}$ and $\varepsilon^q_{j_0\cdots j_{q-k-1}}$ such that
\begin{equation}
\varphi^i\circ\varepsilon^q_{i_0\cdots i_{q-k-1}}=\varphi^j\circ\varepsilon^q_{j_0\cdots j_{q-k-1}},
\end{equation}
then it holds that
\begin{equation}
\Psi^i_s\circ\varepsilon^q_{i_0\cdots i_{q-k-1}}=\Psi^j_s\circ\varepsilon^q_{j_0\cdots j_{q-k-1}}
\quad \mbox{for} \;\forall \; s\in[0,1].
\end{equation}
\end{enumerate}
We call $d_s=\sum n_j\Psi^j_s$ ($s\in[0,1]$) a {\it pre-t-homotopy\/} of $f\circ c$ with respect to $c$. We denote $d$ by the symbol $(f\circ c)^{(c,P)}$.
\end{definition}

Note that the example of Figure 5.1 has no pre-t-modifications.

\begin{lemma}\label{pre-t-lem}
Let $M$ and $N$ be orbifolds, $f:M\rightarrow N$ a b-continuous map, and $c=\sum^{\ell}_{\mu=1}n_{\mu}\varphi^{\mu}$ a $q$-dimensional pre-transversal singular chain of $M$. Let $c'$ be a sum of faces of $c$ such that $f\circ c'$ are already pre-transverse. Then there exists a pre-t-homotopy of $f\circ c$ with respect to $c$ fixing $f\circ c'$.
\end{lemma}

\begin{proof}
We can take atlases $(\{U_{\xi}\}_{\xi\in I},\{\tilde{V}_{\eta}\rightarrow V_{\eta}\}_{\eta\in J},\{(\tilde{\varphi}^{\mu})_{\xi\eta}:U_{\xi}\rightarrow \tilde{V}_{\eta}\}_{\xi\in I,\eta\in J_{\mu\xi}})$ and $(\{\tilde{V}_{\eta}\rightarrow V_{\eta}\}_{\eta\in J},\{\tilde{W}_{\rho}\rightarrow W_{\rho}\}_{\rho\in K},\{\tilde{f}_{\eta\rho}:\tilde{V}_{\eta}\rightarrow \tilde{W}_{\rho}\}_{\eta\in J,\rho\in K_{\eta}})$ of $\varphi^{\mu}$'s and $f$, repectively. Thus $(\{U_{\xi}\}_{\xi\in I},\{\tilde{W}_{\rho}\rightarrow W_{\rho}\}_{\rho\in K},\{\tilde{f}_{\eta\rho}\circ(\tilde{\varphi}^{\mu})_{\xi\eta}:U_{\xi}\rightarrow \tilde{W}_{\rho}\}_{\xi\in I,\rho\in K_{\mu\xi}})$ is an atlas of $f\circ\varphi^{\mu}$, where $K_{\mu\xi}=\{\rho\in K\;|\;|f\circ\varphi^{\mu}|(U_{\xi})\subset W_{\rho}\}$.

To construct a pre-t-homotopy $\Psi_s=\sum^{\ell}_{\mu=1}n_{\mu}\Psi^{\mu}_s$ of $f\circ c$ we will define
\begin{equation}
(\Psi^{\mu}_{(i_0\cdots i_k)})_s:\triangle^q_{(i_0\cdots i_k)}\rightarrow M,\; \mu=1,2,\dots,\ell,\; 0\leq i_0<\cdots<i_k\leq q
\end{equation}
inductively on $k$ such that
\begin{enumerate}
\item[(1)]
\begin{equation*}
\begin{split}
& ((\tilde{\Psi}^{\mu}_{(i_0\cdots i_k)})_0)_{\xi_{(i_0\cdots i_k,0)}\rho_{(i_0\cdots i_k,0)}}|(U_{\xi_{(i_0\cdots i_k,0)}}\cap U_{\xi}) \\
& \qquad =(\tilde{f}_{\eta\rho}\circ(\tilde{\varphi}^{\mu})_{\xi\eta})|(U_{\xi_{(i_0\cdots i_k,0)}}\cap U_{\xi}),
\end{split}
\end{equation*}
\item[(2)] $(\Psi^{\mu}_{(i_0\cdots i_k)})_1$ is pre-transverse,
\item[(3)] $(\Psi^{\mu}_{(i_0\cdots i_k)})_s$'s satisfy the face condition of $(\varphi^{\mu}|\triangle^q_{(i_0\cdots i_k)})$'s,
\item[(4)]
\begin{equation*}
\begin{split}
& ((\tilde{\Psi}^{\mu}_{(i_0\cdots i_k)})_s)_{\xi_{(i_0\cdots i_k,s)}\rho_{(i_0\cdots i_k,s)}}|(U_{\xi_{(i_0\cdots i_k,s)}}\cap U_{(i_0\cdots \check{i}_a\cdots i_k)}) \\
& \qquad =((\tilde{\Psi}^{\mu}_{(i_0\cdots \check{i}_a\cdots i_k)})_s)_{\xi_{(i_0\cdots \check{i}_a\cdots i_k)}\rho_{(i_0\cdots \check{i}_a\cdots i_k)}}|(U_{\xi_{(i_0\cdots i_k,s)}}\cap U_{(i_0\cdots \check{i}_a\cdots i_k)}),
\end{split}
\end{equation*}
\item[(5)] If $(f\circ\varphi^{\mu})|\triangle^q_{(i_0\cdots i_k)}$ is pre-transverse, then
\begin{equation*}
\begin{split}
&((\tilde{\Psi}^{\mu}_{(i_0\cdots i_k)})_s)_{\xi_{(i_0\cdots i_k,s)}\rho_{(i_0\cdots i_k,s)}}|(U_{\xi_{(i_0\cdots i_k,s)}}\cap U_{\xi}) \\
& \qquad =(\tilde{f}_{\eta\rho}\circ(\tilde{\varphi}^{\mu})_{\xi\eta})|(U_{\xi_{(i_0\cdots i_k,s)}}\cap U_{\xi})
\end{split}
\end{equation*}
and \quad $U_{\xi_{(i_0\cdots i_k,s)}}=U_{\xi_{(i_0\cdots i_k,0)}}$,
\end{enumerate}
where
\begin{equation*}
\begin{split}
& (\{U_{\xi_{(i_0\cdots i_k,s)}}\}_{\xi_{(i_0\cdots i_k,s)}\in I_{(i_0\cdots i_k,s)}},
\{W_{\rho_{(i_0\cdots i_k,s)}}\}_{\rho_{(i_0\cdots i_k,s)}\in K_{(i_0\cdots i_k,s)}}, \\
&\quad \{((\tilde{\Psi}^{\mu}_{(i_0\cdots i_k)})_s)_{\xi_{(i_0\cdots i_k,s)}\rho_{(i_0\cdots i_k,s)}} \\
& \quad\qquad :U_{\xi_{(i_0\cdots i_k,s)}}\rightarrow W_{\rho_{(i_0\cdots i_k,s)}}\}_{\xi_{(i_0\cdots i_k,s)}\in I_{(i_0\cdots i_k,s)},\rho_{(i_0\cdots i_k,s)}\in K_{\xi_{(i_0\cdots i_k,s)}}})
\end{split}
\end{equation*}
is an atlas of $(\Psi^{\mu}_{(i_0\cdots i_k)})_s$.

We order $\{(\mu,i_0)\}$ lexicographically and will construct $(\Psi^{\mu}_{i_0})_s$'s, $\mu=1,2,\dots,\ell$, $0\leq i_0\leq q$, in that order as follows:

Let $U_{\xi}$, $V_{\eta}$, and $W_{\rho}$ be local charts of $\triangle^q$, $M$, and $N$, respectively such that $x_{i_0}\in U_{\xi}$, $|\varphi^{\mu}|(U_{\xi})\subset V_{\eta}$, and $|f|(V_{\eta})\subset W_{\rho}$.

1) \quad
If $|f\circ\varphi^{\mu}|(x_{i_0})\not\in\Sigma M$, then we define as
\begin{equation*}
(\tilde{\Psi}^{\mu}_{i_0})_s(x_{i_0}):=(\tilde{f}_{\eta\rho}\circ(\tilde{\varphi}^{\mu})_{\xi\eta})(x_{i_0}).
\end{equation*}

2) \quad
If $|f\circ\varphi^{\mu}|(x_{i_0})\in\Sigma M$, and there exists no $(\nu,j_0)<(\mu,i_0)$ such that $\varphi^{\nu}(x_{j_0})=\varphi^{\mu}(x_{i_0})$.

Let $G_{\rho}$ be the local group of $W_{\rho}$. Since dim {\it Sing} $(G_{\rho})\leq$dim $W_{\rho}-1$, there exists a path in $\tilde{W}_{\rho}$ from $(\tilde{f}_{\eta\rho}\circ(\tilde{\varphi}^{\mu})_{\xi\eta})(x_{i_0})$ to $\tilde{W}_{\rho}-${\it Sing} $(G_{\rho})$. We define $(\tilde{\Psi}^{\mu}_{i_0})_s(x_{i_0})$ by that path.

3) \quad
If $|f\circ\varphi^{\mu}|(x_{i_0})\in\Sigma M$, and for some $(\nu,j_0)<(\mu,i_0)$, $\varphi^{\nu}(x_{j_0})=\varphi^{\mu}(x_{i_0})$. 

Let $U_{\xi '}$ be a local chart of $\triangle^q$ such that $|\varphi^{\nu}|(U_{\xi '})\subset V_{\eta}$ and let $G_{\eta}$ be the local chart of $V_{\eta}$. Since $\varphi^{\nu}(x_{j_0})=\varphi^{\mu}(x_{i_0})$, there exists an element $\sigma\in G_{\eta}$ such that $\sigma(\tilde{\varphi}^{\nu})_{\xi'\eta}(x_{j_0})=(\tilde{\varphi}^{\mu})_{\xi\eta}(x_{i_0})$. From Lemma 3.13(ii), $f$ induce a homomorphism $f_*:G_{\eta}\rightarrow G_{\rho}$. We define
\begin{equation*}
(\tilde{\Psi}^{\mu}_{i_0})_s(x_{i_0}):=f_*(\sigma)(\tilde{\Psi}^{\nu}_{j_0})_s(x_{j_0}).
\end{equation*}

Suppose that we have done for the $(k-1)$-skeleton. Then, by Lemma \ref{sunny1-lem} and the inductive hypothesis (4), we have $((\tilde{\Psi}^{\mu}_{(i_0\cdots i_k)})_s)_{\xi_{(i_0\cdots i_k,s)}\rho_{(i_0\cdots i_k,s)}}$ as the sunny extension of $((\tilde{\Psi}^{\mu}_{(i_0\cdots \check{i}_a\cdots i_k)})_s)_{\xi_{(i_0\cdots \check{i}_a\cdots i_k)}\rho_{(i_0\cdots \check{i}_a\cdots i_k)}}$'s and $(\tilde{f}_{\eta\rho}\circ(\tilde{\varphi}^{\mu})_{\xi\eta})$'s.
\end{proof}

Note that there exist a b-continuous map $f$ and a pre-transversal chain $c$ such that $f\circ c$ has no pre-t-modifiations with respect to \lq\lq $f\circ c$\rq\rq. See Figure 5.4.

\begin{figure}[htbp]
\begin{center}
\includegraphics[width=\linewidth]{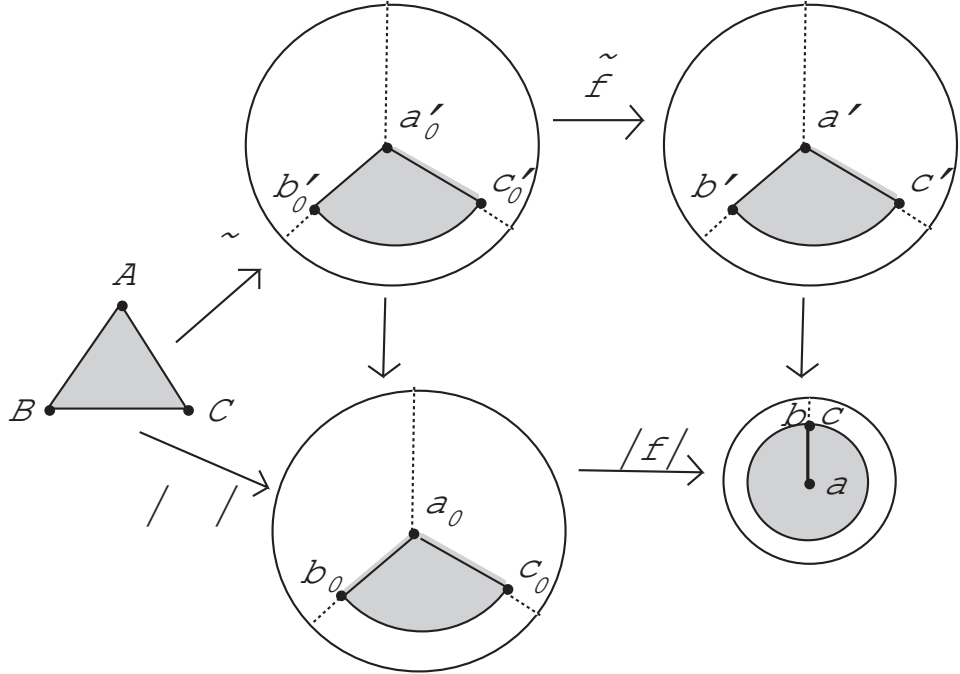}
\end{center}
\caption{$c=\phi=(|\phi|,\tilde{\phi}):\triangle^2\rightarrow D^2,f=(|f|,\tilde{f}):D^2\rightarrow D^2(3)$}
\end{figure}

\begin{lemma}\label{pre-t-partial-lem}
Let $M$ and $N$ be orbifolds, $f:M\rightarrow N$ a b-continuous map, and $c$ a pre-transversal singular chain of $M$. Suppose that $f\circ\partial c$ is pre-transverse.  If $d_s$ is a pre-t-homotopy of $f\circ c$ with respect to $c$ fixing $f\circ\partial c$, then for any $s\in[0,1]$ it holds that $\partial d_s=f\circ\partial c$. In particular, if $c$ is a cycle, then so is $d_s$.
\end{lemma}

\begin{proof}
Similarly to Lemma \ref{tpartial-lem}.
\end{proof}

\begin{theorem}\label{f*def-th}
Let $M$, $N$ be orbifolds and let $f:M\rightarrow N$ be a b-continuous map. Then an induced homomorphism $f_*:t$-$H_q(M)\rightarrow t$-$H_q(N)$ can be defined by $[c]\mapsto[((f\circ c)^{(c,P)})^T]$.
\end{theorem}

\begin{proof}
Take any $c\in t$-$Z_q(M)$.  By Lemma \ref{pre-t-lem} there exists a pre-t-modification $(f\circ c)^{(c,P)}$ of $f\circ c$ with respect to $c$, and by Lemma \ref{t-mod-lem} there exists a t-modification $((f\circ c)^{(c,P)})^T$ of $(f\circ c)^{(c,P)}$. By Lemmas \ref{tpartial-lem} and \ref{pre-t-partial-lem} $((f\circ c)^{(c,P)})^T$ is a cycle.

First take another pre-t-modification $(f\circ c)^{(c,P')}$ of $f\circ c$ with respect to $c$, and any t-modification $((f\circ c)^{(c,P')})^{T'}$ of $(f\circ c)^{(c,P')}$. We will show that the two t-modifications are t-homologue. Put $f\circ c=\sum_{i=1}^{\ell}n_i\varphi^i$ and denote the product of pre-t-homotopy and t-homotopy from $f\circ c$ to $((f\circ c)^{(c,P)})^T$ (resp. $((f\circ c)^{(c,P')})^{T'}$) by $\sum_{i=1}^{\ell}n_i\phi^i_s$ (resp. $\sum_{i=1}^{\ell}n_i\psi^i_s$). Define a continuous map $\Phi^i:\triangle^q\times[-1,1]\rightarrow N$ by
\begin{equation}
\Phi^i(x,s):=
\left\{\begin{array}{ll}
\phi^i_s(x) & \quad (0\leq s\leq 1), \\
\noalign{\vskip2mm}
\psi^i_{-s}(x) & \quad (-1\leq s\leq 0).
\end{array}\right.
\end{equation}
Let $Q^q_k:\triangle^{q+1}\rightarrow \triangle^q\times[0,1]$ be a linear embedding defined by
\begin{equation}
Q^q_k(v_i):=
\left\{\begin{array}{ll}
(v_i,0) & \quad (0\leq i\leq k), \\
(v_{i-1},1) & \quad (k+1\leq i\leq q+1)
\end{array}\right.
\end{equation}
where $v_0,v_1,\dots,v_m$ are the vertices of $\triangle^m$. By taking the composition of each $\Phi^i$ and the sum of $Q^q_k$ we construct a $(q+1)$-dimensional singular chain $d$ of $N$ as
\begin{equation}
d:=\sum_{i=1}^{\ell}n_i\Biggl(\Phi^i\circ\sum_{k=0}^q(-1)^kQ^q_k\Biggr).
\end{equation}
Note that
\begin{equation}\label{partiald-eq}
\partial d=((f\circ c)^{(c,P)})^T-((f\circ c)^{(c,P')})^{T'}.
\end{equation}
Since $\partial d$ is pre-transverse, so is $d$. By Lemma \ref{t-mod-lem} there exists a t-modification $d^T$ of $d$ fixing $\partial d$. By Lemma \ref{tpartial-lem} and \eqref{partiald-eq}
\begin{equation}
\partial (d^T)=\partial d=((f\circ c)^{(c,P)})^T-((f\circ c)^{(c,P')})^{T'}.
\end{equation}

Next take another t-singular cycle $c'$ such that $[c]=[c']$ in $t$-$H_q(M)$. Then there exists a t-singular chain $d'\in t$-$C_{q+1}(M)$ such that $c-c'=\partial d'$. By Lemma \ref{pre-t-lem} there exists a pre-t-modification $(f\circ d')^{(d',P)}$ of $f\circ d'$ with respect to $d'$, and by Lemma \ref{t-mod-lem} there exists a t-modification $((f\circ d')^{(d',P)})^T$ of $(f\circ d')^{(d',P)}$. Then
\begin{equation}
\begin{align}
\partial(((f\circ d')^{(d',P)})^T)
&=(\partial((f\circ d')^{(d',P)}))^T \\
&=((f\circ(\partial d'))^{(\partial d',P)})^T \nonumber \\
&=((f\circ c-f\circ c')^{(c-c',P)})^T \nonumber \\
&=((f\circ c)^{(c,P)})^T-((f\circ c')^{(c',P)})^T. \nonumber
\end{align}
\end{equation}
Thus $f_*$ is well-defined.

To show that $f_*$ is a homomorphism we take any two elements $[a]$, $[b]$ of $t$-$H_q(M)$. Let $(f\circ a)^{(a,P_1)}$ (resp. $(f\circ b)^{(b,P_2)}$) be a pre-t-modification of $f\circ a$ (resp. $f\circ b$), and let $((f\circ a)^{(a,P_1)})^{T_1}$ (resp. $((f\circ b)^{(b,P_2)})^{T_2}$) be a t-modification of $(f\circ a)^{(a,P_1)}$ (resp. $(f\circ b)^{(b,P_2)}$). Since
$((f\circ a)^{(a,P_1)})^{T_1}+((f\circ b)^{(b,P_2)})^{T_2}$ is a t-modification of a pre-t-modification of $f\circ (a+b)$, 
\begin{equation}
\begin{align}
[((f\circ (a+b))^{(a+b,P)})^T]
&=[((f\circ a)^{(a,P_1)})^{T_1}+((f\circ b)^{(b,P_2)})^{T_2}] \\
&=[((f\circ a)^{(a,P_1)})^{T_1}]+[((f\circ b)^{(b,P_2)})^{T_2}]. \nonumber
\end{align}
\end{equation}
\end{proof}

\begin{proposition}\label{f*-prop}
Let $M$, $N$ and $L$ be orbifolds, and let $f:M\rightarrow N$, $g:N\rightarrow L$ be b-continuous maps such that $g\circ f$ is also a b-continuous map. Then $(g\circ f)_*=g_*\circ f_*$.
\end{proposition}

\begin{proof}
Take any $c\in t$-$Z_q(M)$. By Lemmas \ref{t-mod-lem} and \ref{pre-t-lem} there exist products of (pre-)t-homotopies $d_1$ from $g\circ f\circ c$ to $((g\circ f\circ c)^{(c,P_1)})^{T_1}$, $d_2$ from $f\circ c$ to $((f\circ c)^{(c,P_2)})^{T_2}$, and $d_3$ from $g\circ \{((f\circ c)^{(c,P_2)})^{T_2}\}$ to $((g\circ \{((f\circ c)^{(c,P_2)})^{T_2}\})^{((f\circ c)^{(c,P_2)})^{T_2},P_3)})^{T_3}$. Then 
\begin{equation}
\begin{align}
  & \partial(-d_1+g\circ d_2+d_3) \\
=\; & ((g\circ f\circ c)^{(c,P_1)})^{T_1}-((g\circ ((f\circ c)^{(c,P_2)})^{T_2})^{((f\circ c)^{(c,P_2)})^{T_2},P_3)})^{T_3}. \nonumber
\end{align}
\end{equation}
We can make a pre-transversal singular chain $d\in C_{q+1}(L)$ from $-d_1+g\circ d_2+d_3$ such that $\partial d=((g\circ f\circ c)^{(c,P_1)})^{T_1}-((g\circ ((f\circ c)^{(c,P_2)})^{T_2})^{((f\circ c)^{(c,P_2)})^{T_2},P_3)})^{T_3}$. By Lemma \ref{t-mod-lem} there exists a t-modification $d^T$ of $d$. Then by Lemma \ref{tpartial-lem}
\begin{equation}
\partial(d^T)=\partial d=((g\circ f\circ c)^{(c,P_1)})^{T_1}-((g\circ ((f\circ c)^{(c,P_2)})^{T_2})^{((f\circ c)^{(c,P_2)})^{T_2},P_3)})^{T_3}.
\end{equation}
\end{proof}

\begin{theorem}\label{f*-th}
Let $M$, $N$ be orbifolds and let $f,g:M\rightarrow N$ be b-continuous maps. If $f$ and $g$ are homotopic {\rm (}and automatically they are b-homotopic{\rm )}, then $f_*=g_*$.
\end{theorem}

\begin{proof}
Take any $c\in t$-$Z_q(M)$. By Lemmas \ref{t-mod-lem} and \ref{pre-t-lem} there exist products of (pre-)t-homotopies $d_1$ from $f\circ c$ to $((f\circ c)^{(c,P_1)})^{T_1}$, $d_2$ from $g\circ c$ to $((g\circ c)^{(c,P_2)})^{T_2}$. Since $f\circ c$ and $g\circ c$ are homotopic, there exists a prizm from $f\circ c$ to $g\circ c$. Then we can make a pre-transversal singular chain $d\in C_{q+1}(N)$ such that $\partial d=((f\circ c)^{(c,P_1)})^{T_1}-((g\circ c)^{(c,P_2)})^{T_2}$. By Lemma \ref{t-mod-lem} there exists a t-modification $d^T$ of $d$. Then by Lemma \ref{tpartial-lem}
\begin{equation}
\partial(d^T)=\partial d=((f\circ c)^{(c,P_1)})^{T_1}-((g\circ c)^{(c,P_2)})^{T_2}.
\end{equation}
\end{proof}

\begin{theorem}\label{b-homotopy-inv-th}
Let $M$ and $N$ be orbifolds. If they are b-homotopy equivalent, then their t-singular homology groups $t$-$H_q(M)$ and $t$-$H_q(N)$ are isomorphic.
\end{theorem}

\begin{proof}
Since $M$ and $N$ are b-homotopy equivalent, there exist b-continuous maps $f:M\rightarrow N$, $g:N\rightarrow M$ such that $g\circ f$, $f\circ g$ are b-continuous and $g\circ f$ (resp. $f\circ g$) is homotopic to $id_M$ (resp. $id_N$). By Proposition \ref{f*-prop} and Theorem \ref{f*-th} $g_*\circ f_*=(g\circ f)_*=(id_M)_*=id_{t\mbox{-}H_q(M)}$ and $f_*\circ g_*=(f\circ g)_*=(id_N)_*=id_{t\mbox{-}H_q(N)}$.
\end{proof}

\section{Homotopy groups and Hurewicz homomorphisms}

The $q$-th homotopy group $\pi_q(M)$ of an orbifold $M$ was defined in \cite{F-S} as the $(q-1)$-th homotopy group of the (orbifold) loop space of $M$. In this section we regard the $q$-th homotopy group of $M$ as the group of homotopy classes of $q$-dimensional singular spheres with base point $x_0$ by using of b-continuous maps, see Section 3 for the definitions of homotopy and b-continuous map. Then we define the Hurewicz homomorphism from the $q$-th homotopy group $\pi_q(M,x_0)$ to the t-singular homology group $t$-$H_q(M)$, which we defined in Section 4. To define the Hurewicz homomorphism we use some results in Section 5.

\begin{definition}
A $q$-dimensional singular sphere of an orbifold $M$ with base point $x_0\in|M|-\Sigma M$ is a b-continuous map $a=(|a|,\{\tilde{a}_{i\nu}\}):[0,1]^q\rightarrow M$ such that for any $t\in\partial([0,1]^q)$, $|a|(t)=x_0$. A 1-dimensional singular sphere coincides with a loop defined in Definition \ref{loop-def}. The product $a\cdot b$ of $q$-dimensional singular spheres $a,b:[0,1]^q\rightarrow M$ with base point $x_0\in|M|-\Sigma M$ is defined by
\begin{equation}
a\cdot b(t_1,\dots,t_q)=
\left\{
\begin{array}{ll}
a(t_1,\dots,t_{q-1},2t_q)   & 0\leq t_q\leq\frac{1}{\; 2\;}, \\
\noalign{\vskip2mm}
b(t_1,\dots,t_{q-1},2t_q-1) & \frac{1}{\; 2\;}\leq t_q\leq 1.
\end{array}\right.
\end{equation}
\end{definition}

\begin{definition}
Let $M$ be an orbifold and let $x_0$ be a point of $|M|-\Sigma M$. The $q$-dimensional homotopy group $\pi_q(M,x_0)$ of $M$ with base point $x_0$ is the group of homotopy classes, relative to $\partial[0,1]^q$, of $q$-dimensional singular spheres in $M$ with base point $x_0$, where the product of $[a]$ and $[b]$ is defined by
\begin{equation}
[a][b]:=[a\cdot b].
\end{equation}
If $q\geq 2$, $\pi_q(M,x_0)$ is commutative by the definition of products.
\end{definition}

\begin{definition}\label{homotopy-def}
Let $\varphi:[0,1]^q\rightarrow M$ be a $q$-dimensional singular sphere of an orbifold $M$ with base point $x_0$. Let $\mathfrak{S}^q$ be the standard (t-)singular cycle which represents a generator of ($t$-)$H_q(S^q)$. Then we define the Hurewicz homomorphism $\mathfrak{H}_q:\pi_q(M,x_0)\rightarrow t$-$H_q(M)$ by $[\varphi]\mapsto\varphi_*[\mathfrak{S}^q]$ where $\varphi_*$ is defined in Theorem \ref{f*def-th}.
\end{definition}

\begin{theorem}\label{Hure-th}
\begin{enumerate}
\item[(i)] The above map $\mathfrak{H}_q$ is a well-defined homomorphism.
\item[(ii)] If $n=1$ and $M$ is arcwise connected, then the homomorphism $\mathfrak{H}_1$ is surjective and its kernel is the commutator subgroup of $\pi_1(M,x_0)$.
\end{enumerate}
\end{theorem}

\begin{proof}
(i) \quad
By Theorems \ref{f*def-th} and \ref{f*-th} $\mathfrak{H}_q$ is well-defined. To show that $(\varphi\cdot\psi)_*[\mathfrak{S}^q]=\varphi_*[\mathfrak{S}^q]+\psi_*[\mathfrak{S}^q]$ we use the argument parallel to that in the proof of usual Hurewicz homomorphism. Hence $\mathfrak{H}_q$ is a homomorphism.

\noindent
(ii) \quad
It is clear that $\mathfrak{H}_1$ is surjective and that Ker $\mathfrak{H}_1\supset[\pi_1(M,x_0),\pi_1(M,x_0)]$. To show Ker $\mathfrak{H}_1\subset[\pi_1(M,x_0),\pi_1(M,x_0)]$, take any $[\varphi]\in\pi_1(M,x_0)$ such that $[((\varphi\circ\mathfrak{S}^1)^{(\mathfrak{S}^1,P)})^T]=0$ in $t$-$H_1(M)$. By takeing another representative of $[\varphi]$ if necessarily, we may assume that $\varphi\circ\mathfrak{S}^1$ is already transverse with respect to $\Sigma M$. Hence, similarly as in the proof of the usual Hurewicz theorem, we conclude that $[\varphi]\in[\pi_1(M,x_0),\pi_1(M,x_0)]$.
\end{proof}

\section{Exact sequences}

Let $x_0,\dots,x_k$ be the vertices of $\triangle^k$ and let $y_0,\dots,y_k$ be points in $\triangle^n$ which is not necessarily linearly independent. Define a linear map from $\triangle^k$ into $\triangle^n$ by $\sum_{i=0}^k\alpha_ix_i\mapsto\sum_{i=0}^k\alpha_iy_i$ and denote it by $e[y_0\cdots y_k]$. If the points $y_0,\dots,y_k$ are in a general position, this map is an embedding.

\begin{definition}
Let $\Psi:\triangle^q\rightarrow M$, $\phi:\triangle^p\rightarrow M$ be singular simplices of an orbifold $M$ and let $a$ be a point of Int $\triangle^p$. We define the {\it refinement\/} $Sd_{(\phi,a)}(\Psi)$ of $\Psi$ with respect to $(\phi,a)$ as follows:
\begin{enumerate}
\item[(i)] If $\phi$ is the $(i_0\cdots i_p)$-face of $\Psi$, then
\begin{equation}
Sd_{(\phi,a)}(\Psi):=\sum_{k=0}^p(-1)^{i_k+i_p}\Psi\circ e[x_0\cdots\check{x}_{i_k}\cdots x_{i_p}\; a'\; x_{i_p+1}\cdots x_q]
\end{equation}
where $\check{x}_{i_k}$ means deleting $x_{i_k}$ and $a':=e[x_{i_0}\cdots x_{i_p}](a)$.
\item[(ii)] If $\phi$ is not a face of $\Psi$, then $Sd_{(\phi,a)}(\Psi)=\Psi$.
\end{enumerate}
By extending the above definition to $C_q(M)$ linearly we define a {\it refinement operator\/} $Sd_{(\phi,a)}:C_q(M)\rightarrow C_q(M)$ with respect to $(\phi,a)$.
\end{definition}

\begin{definition}
Let $\Psi:\triangle^q\rightarrow M$, $\phi:\triangle^p\rightarrow M$ be singular simplices of an orbifold $M$ and let $a$ be a point of Int $\triangle^p$. We define the {\it prism\/} $P^j_{(\phi,a)}(\Psi)$ {\it of\/} $\Psi$ {\it with respect to\/} $Sd_{(\phi,a)}$ as follows:
\begin{enumerate}
\item[(i)] If $\phi$ is the $(i_0\cdots i_p)$-face of $\Psi$, then, for $0\leq j\leq i_0$,
\begin{equation}
P^j_{(\phi,a)}(\Psi):=\sum_{k=0}^p(-1)^{i_k+i_p}\Psi\circ e[x_0\cdots x_jx_j\cdots x_{i_0}\cdots\check{x}_{i_k}\cdots x_{i_p}\; a'\; x_{i_p+1}\cdots x_q]
\end{equation}
where $\check{x}_{i_k}$ means deleting $x_{i_k}$ and $a':=e[x_{i_0}\cdots x_{i_p}](a)$, and for $i_0+1\leq j\leq n$,
\begin{equation}
P^j_{(\phi,a)}(\Psi):=\Psi\circ e[x_0\cdots x_jx_j\cdots x_q].
\end{equation}
\item[(ii)] If $\phi$ is not a face of $\Psi$, then
\begin{equation}
P^j_{(\phi,a)}(\Psi):=\Psi\circ e[x_0\cdots x_jx_j\cdots x_q].
\end{equation}
\end{enumerate}
Furthermore, we define
\begin{equation}
P_{(\phi,a)}(\Psi):=\sum_{j=0}^q(-1)^{j+1}P^j_{(\phi,a)}(\Psi).
\end{equation}
By extending the above definition to $C_q(M)$ linearly we define a {\it prizm operator\/} $P_{(\phi,a)}:C_q(M)\rightarrow C_{q+1}(M)$ with respect to $Sd_{(\phi,a)}$.
\end{definition}

By a calculation we see the following lemma.

\begin{lemma}\label{prizm-lemma}
Let $\Psi:\triangle^q\rightarrow M$, $\phi:\triangle^p\rightarrow M$ be singular simplices of an orbifold $M$ and let $a$ be a point of\/ {\rm Int} $\triangle^p$. Then the following holds:
\begin{enumerate}
\item[(i)]
$$\partial Sd_{(\phi,a)}(\Psi)=Sd_{(\phi,a)}(\partial\Psi).$$
\item[(ii)]
$$\partial P_{(\phi,a)}(\Psi)=\Psi-Sd_{(\phi,a)}(\Psi)-P_{(\phi,a)}(\partial\Psi).$$
\end{enumerate}
Thus for a singular cycle $c$ of $M$, $Sd_{(\phi,a)}(c)$ and $c$ are homologue.
\end{lemma}

Let $M_1$ and $M_2$ be orbifolds. Note that inclusion maps of orbifolds induce homomorphisms of t-singular chain complices by Remark 4.8. Hence, the general homology theory gives us an exact sequence:

\begin{equation}
\begin{CD}
\cdots @>>> t\mbox{-}H_q(M_1\cap M_2)@> i_* >> t\mbox{-}H_q(M_1)\oplus t\mbox{-}H_q(M_2) \\
@> j_* >> H_q(t\mbox{-}C_q(M_1)+t\mbox{-}C_q(M_2))@> k_* >> \cdots
\end{CD}
\end{equation}
where $i(c)=(c,-c)$, $j(a,b)=a+b$, $k_*[z]=[i^{-1}\partial j^{-1}(z)]$. Furthermore, by using Lemma 7.3, it is proved that if $M_1$ and $M_2$ are suborbifolds of an orbifold $M$ and $\mbox{Int}_M(M_1)\cup \mbox{Int}_M(M_2)=M$, then the natural inclusion map $t$-$C_q(M_1)+t$-$C_q(M_2)\subset t$-$C_q(M_1\cup M_2)$ induces an isomorphism $H_q(t$-$C_q(M_1)+t$-$C_q(M_2))\cong t$-$H_q(M)$. Then we have the Mayer-Vietoris exact sequence:

\begin{theorem}\label{exact-th}
Let $M_1$ and $M_2$ be suborbifolds of an orbifold $M$ and $\mbox{\rm Int}_M(M_1)\cup\mbox{\rm Int}_M(M_2)=M$. Then the following sequence is exact:
\begin{equation}
\begin{CD}
\cdots @>>> t\mbox{-}H_q(M_1\cap M_2)@> i_* >> t\mbox{-}H_q(M_1)\oplus t\mbox{-}H_q(M_2)\\
@> j_* >> t\mbox{-}H_q(M)@> k_* >> t\mbox{-}H_{q-1}(M_1\cap M_2)@>>>\cdots.
\end{CD}
\end{equation}
where $i_*([c])=((i_1)_*[c],-(i_2)_*[c])$, $j_*([a],[b])=(j_1)_*[a]+(j_2)_*[b]$, and $i_{\mu}:M_1\cap M_2\subset M_{\mu}$, $j_{\mu}:M_{\mu}\subset M$ are inclusions. 
\end{theorem}

\section{K\"unneth's Formula}

Let $M$, $N$ be orbifolds and let $p_1:M\times N\rightarrow M$, $p_2:M\times N\rightarrow N$ be the natural projections (see \cite[p.\ 157]{finite} for the product of orbifolds). Observe that for any $c\in t$-$C_k(M\times N)$, $p_1\circ c\in t$-$C_k(M)$ and $p_2\circ c\in t$-$C_k(N)$. Hence we can define a map $\rho_k:t$-$C_k(M\times N)\rightarrow \oplus_{i=0}^k\{(t$-$C_i(M))\otimes(t$-$C_{k-i}(N))\}$ by $\rho_k(c)=\sum_{i=0}^k(p_1\circ c\circ\varepsilon_{i+1\cdots k})\otimes(p_2\circ c\circ\varepsilon_{0\cdots i-1})$, where $\varepsilon_{0,-1}=id$ and $\varepsilon_{k+1,k}=id$. It is proved that $\rho_k$ is a chain homotopy equivalent map by the similar manner as in the usual homology theory. Furthermore, for $\mathbb{Z}$-coefficients general chain complices $C$, $D$, it holds that if $C$ and $H_*(C)$ is free, then $H_*(C)\otimes H_*(D)\cong H_*(C\otimes D)$. Then we obtain the K\"unneth's formula:

\begin{theorem}\label{Kunneth-th}
Let $M$ and $N$ be orbifolds. If either $t$-$H_*(M)$ or $t$-$H_*(N)$ is free, then
\begin{equation}
t\mbox{-}H_k(M\times N)\cong\oplus_{i=0}^k(t\mbox{-}H_i(M)\otimes t\mbox{-}H_{k-i}(N))
\end{equation}
for any $k\geq 0$.
\end{theorem}

\section{Examples}

In this section let us calculate the t-singular homology groups with $\mathbb{Z}$-coefficients of several orbifolds using b-homotopy deformation, exact sequences, and K\"unneth's formula. Note that a discal 2-orbifold $D^2(n)$ with index $n\geq 2$ is not b-homotopy equivalent to a point, see Example 3.7.

\subsection{A lemma for calculations : a discal 2-orbifold case}

\begin{lemma}
For a discal 2-orbifold $D^2(n)$, and the product of $D^2(n)$ and the torus $T^{k-2}$ the following holds:
\begin{equation}
\left\{
\begin{array}{ll}
t\mbox{-}H_2(D^2(n);\mathbb{Z})=0, & \\
t\mbox{-}H_k(D^2(n)\times T^{k-2};\mathbb{Z})=0 & \mbox{\rm for} \quad k\geq 3.
\end{array}
\right.
\end{equation}
\end{lemma}

\begin{proof}
To show $t$-$H_2(D^2(n);\mathbb{Z})=0$, take any $c\in t$-$Z_2(D^2(n);\mathbb{Z})$. Let $\sigma:\triangle^2\rightarrow D^2(n)$ be a b-continuous map such that its structure map $\tilde{\sigma}:\triangle^2\rightarrow D^2$ is a homeomorphism. Let $y$ be the singular point of $D^2(n)$. Denote $c=\sum_{i=1}^ma_i\varphi^i+\sum_{j=1}^rb_j\psi^j$ where $|\varphi^i|(\triangle^2)\not\ni y$ and $|\psi^j|(\triangle^2)\ni y$. Let $f:D^2(n)\rightarrow D^2(n)$ be a b-continuous map which is b-homotopic to $id_{D^2(n)}$, and maps a sufficiently small regular neighbourhood $U$ of the cone point isomorphically onto $D^2(n)$, and maps $D^2(n)-U$ onto $\partial D^2(n)$. By Theorem \ref{f*-th}, $f_*=(id_{D^2(n)})_*$. Hence, by replacing $c$ with $f\circ c$, we may assume that each $\varphi^i(\triangle^2)\subset \partial D^2(n)$ and each $\psi^j=\ell_j\sigma$ in $t$-$C_2(D^2(n);\mathbb{Z})$ for an integer $\ell_j$. Since $\sum_ia_i\varphi^i\in t$-$C_2(\partial D^2(n);\mathbb{Z})$, $[\partial\sum_ia_i\varphi^i]=0$ in $t$-$H_1(\partial D^2(n);\mathbb{Z})$. Hence $[\partial\sum_jb_j\psi^j]=0$ in $t$-$H_1(\partial D^2(n);\mathbb{Z})$. On the other hand $[\partial\sum_jb_j\psi^j]=\sum_jb_j\ell_j[\partial\sigma]$. Since $[\partial\sigma]$ is a nonzero element of $t$-$H_1(\partial D^2(n);\mathbb{Z})\cong\mathbb{Z}$, $\sum_jb_j\ell_j=0$. Thus $c=\sum_ia_i\varphi^i$. Since $t$-$H_2(\partial D^2(n);\mathbb{Z})=0$, there exists a chain $d\in t$-$C_3(\partial D^2(n);\mathbb{Z})$ such that $\partial d=c$. Note that $d\in t$-$C_3(D^2(n);\mathbb{Z})$. Therefore $[c]=0$ in $t$-$H_2(D^2(n);\mathbb{Z})$.

Let $p:D^2\times T^{k-2}\rightarrow D^2(n)\times T^{k-2}$ be the natural covering, and let $[\eta]$ be the fundamental class of $H_k(D^2\times T^{k-2},\partial(D^2\times T^{k-2});\mathbb{Z})$. We can show $t$-$H_k(D^2(n)\times T^{k-2};\mathbb{Z})=0$ for $k\geq 3$ by the similar argument as in the proof of $t$-$H_2(D^2(n);\mathbb{Z})=0$ replacing $\sigma$ by $p\circ \eta$.
\end{proof}

\begin{corollary}
\begin{equation}
t\mbox{-}H_k(D^2(n);\mathbb{Z})=0 \quad \mbox{\rm for} \quad k\geq 2.
\end{equation}
\end{corollary}

\begin{proof}
By Theorem \ref{Kunneth-th} and Lemma 9.1.
\end{proof}

\subsection{A discal 2-orbifold $M:=D^2(n)$}

\begin{equation}
t\mbox{-}H_k(M;\mathbb{Z})\cong
\left\{
\begin{array}{ll}
\mathbb{Z} & \mbox{for} \quad k=0, \\
\mathbb{Z}_n & \mbox{for} \quad k=1, \\
0 & \mbox{for} \quad k\geq 2.
\end{array}
\right.
\end{equation}

Since $|M|-\Sigma M$ is connected, $t$-$H_0(M;\mathbb{Z})\cong \mathbb{Z}$. Take a base point $x_0\in|M|-\Sigma M$. Since $\pi_1(M,x_0)\cong\mathbb{Z}_n$, $t$-$H_1(M;\mathbb{Z})\cong \mathbb{Z}_n$ by Theorem  \ref{Hure-th}. By Corollary 9.2 $t$-$H_k(M;\mathbb{Z})=0$ for $k\geq 2$. By Theorem \ref{b-homotopy-inv-th} we see that $D^2(n)$ is not b-homotopy equivalent to a point, that is, a connected 0-dimensional manifold.

Let $B^3(n)$ be a cyclical ballic 3-orbifold with index $n$. Since $B^3(n)$ is b-homotopy equivalent to $D^2(n)$, their t-homology groups are isomorphic by Theorem \ref{b-homotopy-inv-th}:

\begin{equation}
t\mbox{-}H_k(B^3(n);\mathbb{Z})\cong
\left\{
\begin{array}{ll}
\mathbb{Z} & \mbox{for} \quad k=0, \\
\mathbb{Z}_n & \mbox{for} \quad k=1, \\
0 & \mbox{for} \quad k\geq 2.
\end{array}
\right.
\end{equation}

\subsection{An orientable 2-orbifold $M:=F(m_1,\dots,m_r)$ with genus $g$ and $b$ boundary components}

\begin{equation}
t\mbox{-}H_k(M;\mathbb{Z})\cong
\left\{
\begin{array}{ll}
\mathbb{Z} & \mbox{\rm if} \quad k=0, \\
\underbrace{\mathbb{Z}\oplus\cdots\oplus\mathbb{Z}}_{2g+b-1}\oplus\mathbb{Z}_{m_1}\oplus\cdots\oplus\mathbb{Z}_{m_r} & \mbox{\rm if} \quad k=1, \\
0 & \mbox{\rm if} \quad k=2,\;\partial F\neq\emptyset, \\
\mathbb{Z} & \mbox{\rm if} \quad k=2,\; \partial F=\emptyset, \\
0 & \mbox{\rm if} \quad k\geq 3.
\end{array}
\right.
\end{equation}

If $k=0,1$, the results are derived from the argument similar to that in Subsection 9.2. If $k\geq 2$, the exact sequence yields the results.

\begin{remark}
Let $p_i$ be the singular point of $M$ with index $m_i$, and let $D_i$ be the regular neighbourhood of $p_i$. Let $\sigma_i:\triangle^2\rightarrow D_i$ be a continuous map such that its structure map $\tilde{\sigma}_i:\triangle^2\rightarrow \tilde{D}_i$ is a homeomorphism. Let $\sigma_0$ be a chain in $F-\cup_{i=1}^r\stackrel{\circ}{D}_i$ such that $[\sigma_0]$ generates $H_2(F-\cup_{i=1}^r\stackrel{\circ}{D}_i,\partial(F-\cup_{i=1}^r\stackrel{\circ}{D}_i);\mathbb{Z})$. If $\partial F=\emptyset$, we see that a representative chain of the generator of $t$-$H_2(M;\mathbb{Z})$ is of the form $\ell\sigma_0+\frac{\ell}{m_1}\sigma_1+\cdots+\frac{\ell}{m_r}\sigma_r$, where $\ell$ is the least common multiple of $m_1,\dots,m_r$.
\end{remark}

\subsection{A lemma for calculations : a ballic 3-orbifold case}

\begin{lemma}
Let $M$ be a ballic 3-orbifold $B^3(m_1,m_2,m_3)$. Then the following holds:
\begin{equation}
\left\{
\begin{array}{ll}
t\mbox{-}H_3(M;\mathbb{Z})=0, & \\
t\mbox{-}H_k(M\times T^{k-3};\mathbb{Z})=0 & \mbox{\rm for} \quad k\geq 4.
\end{array}
\right.
\end{equation}
\end{lemma}

\begin{proof}
Similar to Lemma 9.1.
\end{proof}

\begin{corollary}
Let $M$ be a ballic 3-orbifold $B^3(m_1,m_2,m_3)$. Then
\begin{equation}
t\mbox{-}H_k(M;\mathbb{Z})=0 \quad \mbox{\rm for} \quad k\geq 3.
\end{equation}
\end{corollary}

\begin{proof}
By Theorem \ref{Kunneth-th} and Lemma 9.4.
\end{proof}

\subsection{A ballic 3-orbifold $M:=B^3(m_1,m_2,m_3)$}

\begin{equation}
t\mbox{-}H_k(M;\mathbb{Z})\cong
\left\{
\begin{array}{ll}
\mathbb{Z} & \mbox{if} \quad k=0, \\
\mathbb{Z}_2\oplus\mathbb{Z}_2 & \mbox{if} \quad k=1,\;(m_1,m_2,m_3)=(2,2,2n), \\
\mathbb{Z}_2 & \mbox{if} \quad k=1,\;(m_1,m_2,m_3)=(2,2,2n+1), \\
\mathbb{Z}_3 & \mbox{if} \quad k=1,\;(m_1,m_2,m_3)=(2,3,3), \\
\mathbb{Z}_2 & \mbox{if} \quad k=1,\;(m_1,m_2,m_3)=(2,3,4), \\
0 & \mbox{if} \quad k=1,\;(m_1,m_2,m_3)=(2,3,5), \\
0 & \mbox{if} \quad k=2,\;(m_1,m_2,m_3)=(2,2,2n+1), \\
\mathbb{Z}_2 & \mbox{if} \quad k=2,\;(m_1,m_2,m_3)\neq(2,2,2n+1), \\
0 & \mbox{if} \quad k\geq 3.
\end{array}
\right.
\end{equation}

If $k=0,1$, the results are derived from the argument similar to that in Subsection 9.2. Take any $c\in t$-$Z_2(M;\mathbb{Z})$. Note that the image of $c$ does not intersect with the cone point of $M$ by the transversality of $c$. Let $f$ be a continuous map from $M$ onto $M$ which is b-homotopic to $id_M$, and maps a sufficiently small regular neighbourhood of the cone point isomorphically onto $M$. By Theorem \ref{f*-th}, $f_*=(id_M)_*$. Hence, by replacing $c$ with $f\circ c$, we may assume that $c\in t$-$C_2(\partial M;\mathbb{Z})$. By Remark 9.3 $c$ is equivalent to $m(\ell\sigma_0+\frac{\ell}{m_1}\sigma_1+\frac{\ell}{m_2}\sigma_2+\frac{\ell}{m_3}\sigma_3)$ for an integer $m$. On the other hand, it is clear that if $c\in t$-$B_2(M)$, then $c$ is equivalent to $m'(n_0\sigma_0+\frac{n_0}{m_1}\sigma_1+\frac{n_0}{m_2}\sigma_2+\frac{n_0}{m_3}\sigma_3)$, where $m'$ is an integer and $n_0$ is the index of the cone point $\Sigma^{(0)} M$. If $(m_1,m_2,m_3)=(2,2,2n+1)$, $n_0=\ell=4n+2$. Then $t$-$H_2(M;\mathbb{Z})=0$. If $(m_1,m_2,m_3)\neq(2,2,2n+1)$, $n_0=2\ell$. Then $t$-$H_2(M;\mathbb{Z})\cong\mathbb{Z}_2$.

By Corollary 9.5 $t$-$H_k(M;\mathbb{Z})=0$ for $k\geq 3$.

\section{singular homology and s-singular homology}

In this section we give a proof of Theorem \ref{under-H}, though it is intuitively clear, by using of a modification and a refinement of a singular chain.

\begin{definition}\label{s-sing-def}
A $q$-dimensional {\it s-singular simplex\/} of an orbifold $M$ is a $q$-dimensional singular simplex $\Psi:\triangle^q\rightarrow M$ of $M$ which satisfies the following:
\begin{enumerate}
\item[(a)] For each $k$, there exists a face $\sigma$ of $\triangle^q$ (possibly, $\triangle^q$ or $\emptyset$) such that $|\Psi|^{-1}(\cup_{j=0}^k\Sigma^{(j)} M)=\sigma$.
\end{enumerate}
The following is easily derived from (a):
\begin{enumerate}
\item[(i)] For each face $\eta$ of $\triangle^q$ there exists a stratum $Q$ of $M$ (possibly, $Q\not\subset\Sigma$) such that $|\Psi|($Int $\eta)\subset Q$.
\item[(ii)] For each $k$ there exists at most one stratum of $k$-dimension that intersects with $|\Psi|(\triangle^q)$.
\end{enumerate}
A $q$-dimensional {\it s-singular chain\/} of $M$ is a finite linear combination $\sum_jn_j\Psi^j$ of s-singular simplices $\Psi^j$ of $M$. We denote the free abelian group with basis of all $q$-dimensional s-singular simplices of $M$ by $s$-$C_q(M)$. By the above definition of s-singular simplex it holds that $\partial(s$-$C_q(M))\subset s$-$C_{q-1}(M)$, and $s$-$C_*(M)$ is a subchain complex of $C_*(M)$. We define the {\it s-singular homology group\/} of $M$ by $s$-$Z_q(M)/s$-$B_q(M)$ and denote it by $s$-$H_q(M)$ where
\begin{equation}
\begin{array}{ll}
s\mbox{-}Z_q(M)&=\{c\in s\mbox{-}C_q(M)\;|\; \partial(c)=0 \} \\
s\mbox{-}B_q(M)&=\{ c\in s\mbox{-}C_q(M) \;|\; \exists d\in s\mbox{-}C_{q+1}(M) \;\mbox{s.t.}\; \partial(d)=c \}.
\end{array}
\end{equation}
\end{definition}

For any singular chain $c=\sum n_i \Psi^i$ of an orbifold $M$ we denote the singular chain $\sum n_i|\Psi^i|$ of $|M|$ by $|c|$ where $|\Psi^i|$ is the underlying map of $\Psi^i$. 

\begin{lemma}\label{refine-lem}
Let $M$ be an orbifold and let $v=\sum_i m_i \varphi^i$ be a $q$-dimensional singular chain of $|M|$ where each $\varphi^i$ is a singular simplex of $|M|$ satisfying {\rm (a)} of Definition \ref{s-sing-def}. Then there exists an s-singular chain $c$ of $M$ with $|c|=v$.
\end{lemma}

\begin{proof}
All we have to do is to construct a continuous map $\Psi^i:\triangle^q\rightarrow M$ such that $|\Psi^i|=\varphi^i$ for each $\varphi^i:\triangle^q\rightarrow |M|$. Let $\{\tilde{V}_{\nu}\rightarrow V_{\nu}\}$ (resp. $\{\tilde{U}_{\mu}\rightarrow U_{\mu}\}$) be a system of local charts of $\triangle^q$ (resp. $M$). Take any $V_j\in\{V_{\nu}\}$ and any $U_k\in\{U_{\mu}\}$ with $\varphi^i(V_j)\subset U_k$. By the assumption we can take a lift $(\tilde{\Psi}^i)_{jk}$ of $\varphi^i|V_j\rightarrow U_k$ to $\tilde{U}_k$. Put $\Psi^i=(\varphi^i,\{(\tilde{\Psi}^i)_{jk}\})$. Then $c=\sum_i m_i \Psi^i$ is a desired chain.
\end{proof}

\begin{proposition}\label{underref-prop}
Let $c$ be a singular chain of $|M|$. Let $c_1$ be a sum of faces of simplices in $c$, each of which satisfies {\rm (a)} of Definition \ref{s-sing-def}. Then there exists a modification $c'$ obtaind from $c$ by homotopies and refinements, which satisfies {\rm (a)} of Definition \ref{s-sing-def} and fixes $c_1$.
\end{proposition}

\begin{definition}
Let $c$ be a singular chain of an orbifold $M$. A singular chain $c'$ of $M$ is called an {\it s-modification\/} of $c$ if $c'$ is obtained from $c$ by homotopies and refinements and is an s-singular chain of $M$. Note that if $c$ and $c'$ are cycles, then they are homologue as singular chains of $M$ by Lemma \ref{prizm-lemma}.
\end{definition}

\begin{proposition}\label{refine-prop}
Let $c$ be a singular chain of an orbifold $M$. Let $c_1$ be a sum of faces of simplices in $c$, each of which is an s-simplices. Then there exists an s-modification of $c$ fixing $c_1$.
\end{proposition}

The following theorem says that the (s-)singular homology of an orbifold does not respect the orbifold structure. On the other hand the t-singular homology does. See Section 9, where the t-singular homology is illustrated through simple but important examples.

\begin{theorem}\label{under-H}
Let $M$ be an orbifold. For each $q$ the following holds:
\begin{equation}
s\mbox{-}H_q(M)\cong H_q(M)\cong H_q(|M|).
\end{equation}
\end{theorem}

\begin{proof}
We define a map $\mu:H_q(M)\rightarrow H_q(|M|)$ by $\mu([c]):=[|c|]$. It is well-defined. Indeed, if $[c]=[c']$, there exists a $(q+1)$-dimensional chain $d\in C_{q+1}(M)$ such that $c'=c+\partial d$. Then
\begin{equation}
\mu([c'])
=\mu([c+\partial d])
=[|c+\partial d|]
=[|c|+|\partial d|]
=[|c|+\partial|d|]
=[|c|].
\end{equation}

Also the map $\mu$ is a homomorphism since
\begin{equation}
\mu([c_1+c_2])
=[|c_1+c_2|]
=[|c_1|+|c_2|]
=[|c_1|]+[|c_2|]
=\mu([c_1])+\mu([c_2]).
\end{equation}

To show that $\mu$ is injective, let $\mu([c])=0$. Take an s-modification $c'$ of $c$. Since $\mu[c]=\mu[c']=[|c'|]=0$, there exists a singular chain $d\in C_{q+1}(|M|)$ such that $|c'|=\partial d$. By Proposition \ref{underref-prop} we can take a modification $d'$ of $d$ which satisfies (a) of Definition \ref{s-sing-def} and fixes $\partial d$. By Lemma \ref{refine-lem} there exists a singular chain $f\in C_{q+1}(M)$ such that $|f|=d'$ and $\partial f=c'$. Thus $[c]=[c']=0\in C_q(M)$.

To show that $\mu$ is surjective, take any element $d\in C_q(|M|)$. By Proposition \ref{underref-prop} there exists a modification $d'$ of $d$ which satisfies (a) of Definition \ref{s-sing-def}. By Lemma \ref{refine-lem} there exists a singular chian $c\in C_q(M)$ such that $|c|=d'$. Thus
\begin{equation}
\mu([c])=[|c|]=[d']=[d].
\end{equation}
Similarly we can show that $s$-$H_q(M)\cong H_q(|M|)$.
\end{proof}

\section{rational coefficients}

In this section we show that the t-singular homology group with rational coefficients of an orbifold is isomorphic to the usual singular homology group with rational coefficients of its underlying space.

\begin{definition}
Let $c=\sum_{i=1}^{\ell}n_i\varphi^i\in C_q(M)$ be a singular chain of an orbifold $M$. Put $\partial c=\sum_{i=1}^rm_i\sigma^i$ ($m_i\neq 0$). Let $\varphi^i_s$, $i=1,2,\dots,\ell$ be homotopies with $\varphi^i_0=\varphi^i$. We call $c_s=\sum_{i=1}^rn_i\varphi^i_s$ a homotopy of $c$ preserving the boundary condition if
\begin{equation}
\partial(c_s)=\sum_{i=1}^rm_i\sigma^i_s
\end{equation}
where $\sigma^i_s=\varphi^j_s\circ\varepsilon_k$ if $\sigma^i=\varphi^j\circ\varepsilon_k$, $j=1,2,\dots,\ell$, $k=0,1,\dots,q$ for each $i=1,2,\dots,r$.
\end{definition}

\begin{lemma}\label{w-t-mod-lem}
Let $M$ be an orbifold and let $c=\sum_{i=1}^k\alpha_i\varphi^i$ be a $q$-dimensional singular chain of $M$. If for each vertex $x_{\nu}$ of $\triangle^q$, and for each $i=1,2,\dots,k$, $\alpha_i$ is a multiple of $w(\varphi^i(x_{\nu}))$, then we can construct a homotopy $c_s$ of $c$ preserving the boundary condition such that $c_1$ is a pre-transversal chain $c^{mP}$.
\end{lemma}

\begin{proof}
Put $y=\varphi^1(x_0)$. We may assume that each of $\varphi^1,\dots,\varphi^{\ell}$ maps at least one vertex to $y$, and none of $\varphi^{\ell+1},\dots,\varphi^k$ maps any vertices to $y$. Take a reduced chart $\tilde{U}\rightarrow U$ around $y$ and let $\tilde{y}$ be the point in $\tilde{U}$ corresponding to $y$. Let $\tilde{\varphi}^i_1,\dots,\tilde{\varphi}^i_{r_i}$, $i=1,2,\dots,\ell$, be the structure maps of $\varphi^i$ each of which maps at least one vertex to $\tilde{y}$, and let $\tilde{\varphi}^i_{r_i+1},\dots,\tilde{\varphi}^i_{s_i}$, $i=1,2,\dots,\ell$, be the others. By the hypothesis we may assume that the structure maps of $\alpha_i\varphi^i$ is of the form
\begin{equation}
\{\alpha'_ig_j\tilde{\varphi}^i_{\mu}\}_{\mu=1,\dots,r_i,g_j\in G_{\tilde{y}}}
\cup\{\alpha_i\tilde{\varphi}^i_{\eta}\}_{\eta=r_i+1,\dots,s_i}
\end{equation}
where $\alpha_i'=\alpha_i/\# G_{\tilde{y}}$. Thus we can slightly modify each $g_j\tilde{\varphi}^i_{\mu}$, $g_j\in G_{\tilde{y}}$, $\mu=1,2,\dots,r_i$, $i=1,2,\dots,\ell$, in a small neighbourhood of $\tilde{y}$ preserving the face conditions of $g_j\tilde{\varphi}^i_{\mu}$'s (not of $\varphi^i$'s) so that it does not map any vertices to $\tilde{y}$ (see Figure 11.1). Then we change the chain $c$ to one whose structure maps are the result of the modification. Modifying on the other vertices as well, we obtain the desired homotopy $c_s$.
\end{proof}

\begin{figure}[htbp]
\begin{center}
\includegraphics[width=\linewidth]{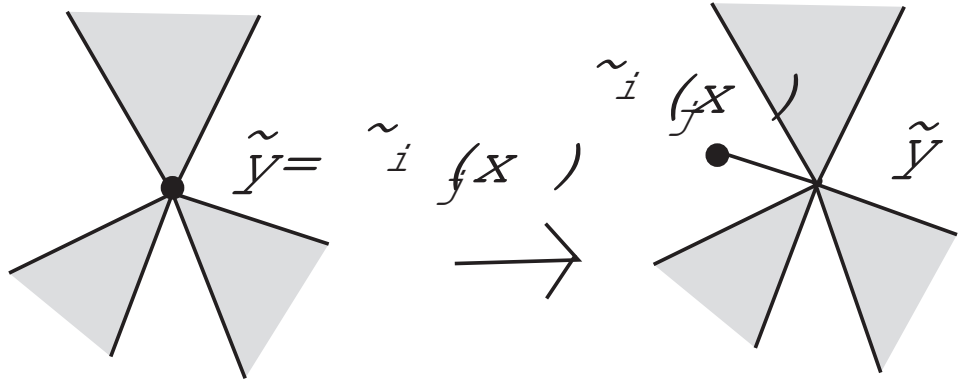}
\end{center}
\caption{a modification of the structure maps in a small neighbourhood of $\tilde{y}$, $i=1,2,\dots,6$, $j\in\{ 0,1,2\}$}
\end{figure}

We call the above $c^{mP}$ a {\it pre-t-modification of $c$ from the multiplicity}.

\begin{remark}
In Lemma 11.2, $c^{mP}$ satisfies the following:
\begin{enumerate}
\item[(i)] If $\partial c$ is pre-transverse, then $\partial c=\partial(c^{mP})$.
\item[(ii)] If $c$ is a cycle, then $[c]=[c^{mP}]$ in $H_q(M)$.
\item[(iii)] By (11.1) we have $\partial(c_s)=(\partial c)_s$ where $(\partial c)_s$ is a homotopy of $\partial c$ and is the restriction of $c_s$.
\end{enumerate}
\end{remark}

\begin{theorem}
Let $M$ be an orbifold. Then for any $q$ the following homology groups with rational coefficients are isomorphic:
\begin{equation}
s\mbox{-}H_q(M;\mathbb{Q})\cong H_q(M;\mathbb{Q})\cong t\mbox{-}H_q(M;\mathbb{Q})\cong H_q(|M|;\mathbb{Q}).
\end{equation}
\end{theorem}

\begin{proof}
For any chain $\sigma=\sum_{i=1}^{\ell}n_i\psi^i\in C_r(M;\mathbb{Q})$, let $w_{\sigma}$ be the least number such that each $n_iw_{\sigma}$ is a multiple of $w(\psi^i(x_{\nu}))$ for $i=1,\dots,\ell$ and $\nu=0,1,\dots,r$ where $x_{\nu}$ are the vertices of $\triangle^r$. By Theorem 10.6, $s$-$H_q(M;\mathbb{Q})\cong H_q(M;\mathbb{Q})\cong H_q(|M|;\mathbb{Q})$. We will show that $H_q(M;\mathbb{Q})\cong t$-$H_q(M;\mathbb{Q})$. Take any $[c]\in H_q(M;\mathbb{Q})$ and put $c=\sum_{i=1}^k\alpha_i\varphi^i\in Z_q(M;\mathbb{Q})$. Since
\begin{equation}
w_cc=\sum_{i=1}^k\alpha_iw_c\varphi^i,
\end{equation}
there exists a pre-t-modification $(w_cc)^{mP}$, which is a cycle by Lemma \ref{w-t-mod-lem}. By Lemma \ref{t-mod-lem} there exists a t-modifiation $((w_cc)^{mP})^T$, which is a cycle. We define a map $\Psi:H_q(M;\mathbb{Q})\rightarrow t$-$H_q(M;\mathbb{Q})$ by
\begin{equation}
\Psi([c])=[\frac{1}{\; w_c\;}((w_cc)^{mP})^T].
\end{equation}

To show that the map $\Psi$ is well-defined, first take a t-modification $((w_cc)^{mP'})^{T'}$ of another pre-t-modification of $w_cc$. By the argument similar as in the proof of Theorem \ref{f*def-th} we can show that $[\frac{1}{\; w_c\;}((w_cc)^{mP})^T]=[\frac{1}{\; w_c\;}((w_cc)^{mP'})^{T'}]\in t$-$H_q(M;\mathbb{Q})$.

Next take any $c,c'\in Z_q(M;\mathbb{Q})$ such that $[c]=[c']\in H_q(M;\mathbb{Q})$. Then there exists $d\in C_{q+1}(M;\mathbb{Q})$ such that $c-c'=\partial d$. Let $\beta$ be the least number such that $\beta w_d$ is a common multiple of $w_c$ and $w_{c'}$. By Lemma \ref{w-t-mod-lem} there exists a pre-t-modification $(\beta w_dd)^{mP}$ of $\beta w_dd$ from the multiplicity, and by Lemma \ref{t-mod-lem} there exists a t-modification $((\beta w_dd)^{mP})^T$ of $(\beta w_dd)^{mP}$. Then, by (iii) of Remark 11.3
\begin{equation}
\begin{align}
\partial(((\beta w_dd)^{mP})^T)
&=(\partial((\beta w_dd)^{mP}))^T \\
&=((\partial(\beta w_dd))^{mP})^T \nonumber \\
&=((\beta w_dc)^{mP}-(\beta w_dc')^{mP})^T \nonumber \\
&=((\beta w_dc)^{mP})^T-((\beta w_dc')^{mP})^T, \nonumber
\end{align}
\end{equation}
which shows
\begin{equation}
[((\beta w_dc)^{mP})^T]=[((\beta w_dc')^{mP})^T].
\end{equation}
Since $w_d$ is a common multiple of $w_c$ and $w_{c'}$,
\begin{equation}
[\frac{\beta w_d}{\; w_c\;}((w_cc)^{mP})^{T}]
=[\frac{\beta w_d}{\; w_{c'}\;}((w_{c'}c')^{mP})^{T}].
\end{equation}
Then we have $[\frac{1}{\; w_c\;}((w_cc)^{mP})^{T}]=[\frac{1}{\; w_{c'}\;}((w_{c'}c')^{mP})^{T}]$. Thus $\Psi$ is well-defined.

To show that $\Psi$ is a homomorphism, take any $[a],[b]\in H_q(M;\mathbb{Q})$. Let $\gamma$ be the least number such that $\gamma w_aw_b$ is a multiple of $w_{a+b}$.
\begin{equation}
\begin{align}
\Psi([a])+\Psi([b])
&=[\frac{1}{\; w_a\;}((w_aa)^{mP_1})^{T_1}]+[\frac{1}{\; w_b\;}((w_bb)^{mP_2})^{T_2}] \\
&=\frac{1}{\;\gamma w_aw_b\;}[((\gamma w_aw_ba)^{mP_1})^{T_1}]+\frac{1}{\;\gamma  w_aw_b\;}[((\gamma w_aw_bb)^{mP_2})^{T_2}] \nonumber \\
&=\frac{1}{\;\gamma  w_aw_b\;}[((\gamma w_aw_ba)^{mP_1})^{T_1}+((\gamma w_aw_bb)^{mP_2})^{T_2}] \nonumber \\
&=\frac{1}{\;\gamma  w_aw_b\;}[((\gamma w_aw_b(a+b))^{mP})^T] \nonumber \\
&=\frac{1}{\;\gamma  w_aw_b\;}[\frac{\gamma w_aw_b}{w_{a+b}}((w_{a+b}(a+b))^{mP})^T] \nonumber \\
&=[\frac{1}{\; w_{a+b}\;}((w_{a+b}(a+b))^{mP})^T] \nonumber \\
&=\Psi([a+b]) \nonumber\\
&=\Psi([a]+[b]). \nonumber
\end{align}
\end{equation}

To show that $\Psi$ is injective take any $[c]\in H_q(M;\mathbb{Q})$ such that $\frac{1}{\; w_c\;}[((w_cc)^{mP})^T]=0\in t$-$H_q(M;\mathbb{Q})$. Then there exists 
$d\in t$-$C_{q+1}(M;\mathbb{Q})$ such that $\partial d=((w_cc)^{mP})^T$. By (ii) of Remark 11.3,  $[w_cc]=[(w_cc)^{mP}]$ in $H_q(M;\mathbb{Q})$. By making a prizm from a t-homotopy, we see that $[(w_cc)^{mP}]=[((w_cc)^{mP})^T]$ in $H_q(M;\mathbb{Q})$. Since $d\in C_{q+1}(M;\mathbb{Q})$, $[w_cc]=0\in H_q(M;\mathbb{Q})$. Thus $[c]=0\in H_q(M;\mathbb{Q})$.

To show that $\Psi$ is surjective, take any $[c]\in t$-$H_q(M;\mathbb{Q})$. Since $c\in Z_q(M;\mathbb{Q})$, $\Psi([c])=[c]$.
\end{proof}

\section{ws-singular cohomology}

Let $M$ be an orbifold and let $\Psi:\triangle^q\rightarrow M$ be an s-singular simplex of $M$. We define the {\it weight\/} $w(\Psi)$ of $\Psi$ as
\begin{equation}
w(\Psi)=w(\Psi(p))
\end{equation}
where $p$ is an interior point of $\triangle^q$. It is well-defined by Definition \ref{s-sing-def}. We define ws-singular cochain $ws$-$C^q(M;\mathbb{Z})$ as
\begin{equation}
ws\mbox{-}C^q(M;\mathbb{Z})=\{\eta\in s\mbox{-}C^q(M;\mathbb{Z})\;|\; \langle\eta,\Psi\rangle\in w(\Psi)\mathbb{Z}\; \mbox{for each $q$-dimensional s-simplex}\; \Psi\}.
\end{equation}
Let $\delta$ be the coboundary operator $\delta:s$-$C^q(M)\rightarrow s$-$C^{q+1}(M)$ defined by
\begin{equation}
\langle\delta(\eta),c\rangle=\langle\eta,\partial c\rangle \quad \mbox{for}\;\eta\in s\mbox{-}C^q(M),c\in s\mbox{-}C_{q+1}(M).
\end{equation}
By restricting $\delta$ to $ws$-$C^q(M)$, we have the cochain complex $(\{ ws$-$C^q(M)\},\delta)$. We define the ws-singular cohomology group $ws$-$H^q(M;\mathbb{Z})$ as the cohomology group of $(\{ ws$-$C^q(M)\},\delta)$. In \cite{ws-coh}, for a compact, orientable $n$-dimensional orbifold $M$, we will give the proof of the Poincar\'e duality theorem:

\begin{equation}
\begin{align}
t\mbox{-}H_{n-q}(M,\partial M;\mathbb{Z})&\cong ws\mbox{-}H^q(M;\mathbb{Z}), \\
ws\mbox{-}H^{n-q}(M,\partial M;\mathbb{Z})&\cong t\mbox{-}H_q(M;\mathbb{Z}).
\end{align}
\end{equation}

\end{document}